\documentclass[11pt]{article}
\usepackage{graphicx}
\usepackage{latexsym,amssymb,amsmath,amscd,amsfonts,mathabx}
\topmargin - 2cm
\oddsidemargin - 0.1cm
\evensidemargin - 5cm
\newtheorem{theorem}{Theorem}[section]

\newtheorem{lemma}[theorem]{Lemma}
\newtheorem{cor}[theorem]{Corollary}

\newtheorem{remark}[theorem]{Remark}
\newtheorem{proof}[theorem]{Proof}

\newcommand{\eproof}{\rule{0.2cm}{0.2cm}}

\newcommand{\re}{\mathbb{R}}

\begin{document}

\title{\bf On systems of fractional nonlinear partial differential equations}

\author{Ravshan Ashurov$^{1,2}$,  \  Oqila Mukhiddinova$^{1,3}$}
\date{}
\maketitle

\begin{center}
{\it $^1$V.I. Romanovskiy Institute of Mathematics, Uzbekistan Academy of Science, University str.,9, Olmazor district, Tashkent, 100174, Uzbekistan; e-mail ashurovr@gmail.com\\
$^2$ Department of Mathematics, New Uzbekistan University, Movarounnahr Street 1, Tashkent, 100000, Uzbekistan; e-mail ashurovr@gmail.com \\
$^3$ Tashkent University of Information Technologies named after Muhammad al-Khwarizmi,
108 Amir Temur Avenue,Tashkent, 100200, Uzbekistan; e-mail oqila1992@mail.ru}

\end{center}
\vspace{10pt}

\begin{abstract}
The work considers a system of fractional order partial differential equations. The existence and uniqueness theorems for the classical solution of initial-boundary value problems are proved in two cases: 1) the right-hand side of the equation does not depend on the solution of the problem and 2)  it depends on the solution, but at the same time satisfies the classical Lipschitz condition with respect to this variable and an additional condition which guarantees a global existence of the solution. Sufficient conditions are found (in some cases they are necessary) on the initial function and on the right-hand side of the equation, which ensure the existence of a classical solution.   In previously known works, linear but more general systems of fractional pseudodifferential equations were considered and the existence of a weak solution was proven in the special classes of distributions.

\end{abstract}

%
\vspace{2pc}

%
%
%
%

\section{Introduction}\label{sec:1}

Over the past two decades, specialists have intensively studied systems of fractional differential equations (ordinary and in partial derivatives) and discovered rich applications, for example, in modeling of processes in biosystems \cite{DasGupta}, \cite{Rihan}, \cite{GuoFang} ecology \cite{Khan}, \cite{Rana}, epidemiology \cite{Zeb}, \cite{Islam}, \cite{Lyapunov-2}, \cite{Goufo}, \cite{Almeida}, \cite{Rajagopal}, etc.

In the study of the properties of the solution and in the approximate calculation of the solution of fractional order ordinary or partial
differential equations representation formulas for the solution play a key role. That is why the research of many authors is devoted to finding such a formula for various systems of equations.
In the case of systems of ordinary differential equations of fractional order, obtaining representation formulas is comparatively simpler. 

One of the first works where a representation of the solution was obtained is the article \cite{Varsha} of the authors Varsha Daftardar-Gejji  and A. Babakhani, where the authors examined the system
\[
D_t^{\alpha}[u(t)-u(0)]=Au(t), \quad u(0)=u_0,
\] 
of time-fractional ordinary differential equations. Here $u(t)$ is a vector-function, $A$ is a nonsingular matrix and $\alpha \in (0,1)$ is scalar. First, the authors considered the case when $A$ has real and distinct eigenvalues. Then the authors moved on to an arbitrary matrix and, using the Jordan form of the matrix, proved the existence and uniqueness of the solution in a small neighborhood in time. A representation for solving the Cauchy problem in the form $u(t)=E_{\alpha}(t^{\alpha}A)u_0$ is obtained, where $E_{\alpha}(\mathcal{Z})$ is the matrix-valued Mittag-Leffler function of a matrix $\mathcal{Z}.$ The authors also investigated the nonlinear system by replacing $Au(t)$ with $f(t, u)$. Reducing the corresponding Cauchy problem to the Volterra integral equation, it is shown that if $f(t, u)$ satisfies the standard Lipschitz condition for the variable $u$, then the integral equation has a unique solution on a sufficiently small interval. Note that in this work the equivalence of the Cauchy problem with the Voltaire equation was not proven. This statement was proven in work \cite{BonillaKilbas} for a more general nonlinear system of equations, where the existence and uniqueness of a generalized solution was shown.

In works \cite{DengLiGuo}, \cite{Odabat} two classes of linear fractional differential systems $D_t^{\beta} u(t) = A u(t)$ are considered, where $A$ is a $m\times m$ matrix with constant coefficients and the orders are equal to a real number (commensurate fractional order) or rational numbers (incommensurate fractional order) between zero and one. In both papers asymptotic stability theorems for the solutions of the Cauchy problem are proved. The authors of work \cite{DengLiGuo} proved the asymptotic stability of the solution also for systems of differential equations with a nonlinear right-hand side. We emphasize that in this work a representation for solving the problem is not written out. Whereas in the work \cite{Odabat} formulas were also obtained to represent the solution to the problem under consideration.

The authors of both papers in the case when the components $\beta_j$ of the vector-order $\beta=(\beta_1, \cdots, \beta_m)$ of the equation are rational, reduced the corresponding system, as in the case of systems with integer order derivatives, to a system with the same fractional order in each equation.

The work \cite{Mamchuev} of M. O. Mamchuev examines a boundary-value problem for a system of multi-term differential equations of fractional order
\[
\sum\limits_{j=1}^n A_j D^{\alpha_j}_{x_j} u(x)=B u(x)+ f(x),
\]
where $A$ and $B$ are $m\times m$ matrices with constant coefficients, $u$ and $f$ are vector functions of $x\in \Omega$ and $\Omega$ is a $n$ dimensional  rectangular domain. The author constructed the Green's function using the Wright function.

Let us note one more work \cite{NewMethod} (see also  various numerical methods in \cite{Erturk}, \cite{Wang}, \cite{Abdulaziz}, \cite{Lyapunov-1}, \cite{Wang}), where the authors propose a method for solving systems of linear and nonlinear fractional partial differential equations that is effective for numerical solution. This method is a combination of the Laplace transform and iterative method. The method does not contain rounding errors, which leads to a reduction in numerical calculations. The paper also presents numerical examples illustrating the accuracy and efficiency of this method.

In the present work, using the method developed in the papers \cite{UAChen, UmarovToAppear, UmarovNew}, the existence and uniqueness of a classical solution to systems of fractional partial differential equations with linear and nonlinear right-hand sides is proven. Namely, we will be dealing with the following system of nonlinear partial differential equations of fractional order

    \begin{equation} \label{system_main}
     \begin{cases} 
	 D_t^{\beta}u_1(t,x) + A_{1, 1}(D) u_1(t,x) +\dots +A_{1, m}(D) u_m(t,x) = h_1(t,x, u_1, \cdots, u_m), \\ 
         D_t^{\beta}u_2 (t,x)+ A_{2, 1}(D) u_1 (t,x)+\dots +A_{2, m} (D)u_m(t,x) = h_2(t,x, u_1, \cdots, u_m), \\ 
          {\cdots } \\
         D_t^{\beta}u_m (t,x) + A_{m, 1}(D) u_1(t,x) +\dots +A_{m, m}(D) u_m (t,x) = h_m(t,x, u_1, \cdots, u_m),     
           \end{cases}
     \end{equation}
where $(t, x)\in \mathbb{R}_+\times \mathbb{R}^n$, $D_t^{\beta}$ is the fractional order derivative of order $ 0<\beta \le 1$ in the sense of Caputo, and $A_{j,k}(D)$ are differential operators with constant coefficients  (described later in Section 2.3). Functions $h_j$ generally speaking, depend nonlinearly on $\langle{u_1, u_2, \dots, u_m \rangle}$; we give exact conditions for these functions in Section 4. 

Linear systems of fractional partial differential equations have been studied by many authors; a review of published work up to 2019 can be found in papers of Anatoly N. Kochubei \cite{Kochuei1} (fractional-parabolic equations, see also \cite{Kochubei3} and \cite{Kochubei4}) and \cite{Kochubei2} (fractional-hyperbolic equations). In  \cite{Kochuei1}, in particular, the Green's functions of the Cauchy problem for linear systems were constructed and their estimates were obtained when the differential expressions $A_{j,k}(D)$ depend on the spatial variables $x$. In all the works listed in \cite{Kochuei1}, the order of all operators $A_{j,k}(D)$ in the system  are the same.

A system of fractional pseudodifferential equations of the form (\ref{system_main}) (when matrix $\{A_{j,k}\}$ consists of pseudodifferential operators) was systematically studied by S. Umarov (see \cite{UmarovToAppear} and \cite{UmarovNew}). In work \cite{UmarovToAppear}, the author obtained representations of the solution to the Cauchy problem and proved the existence of a weak solution in the special classes of distributions in the following three cases: 1) when $\beta\in (0,1]$ takes the same values in all equations of the system, 2) when $\beta$ is a vector $\beta=\langle{\beta_1, \cdots, \beta_m\rangle}$ and $\beta_j\in (0, 1]$ takes different rational values, 3) $\beta_j\in (0, 1]$ takes different real values, but the matrix $\{A_{j,k}\}$ is triangular. In this work, the author also investigated the question of the existence of a strong solution from Sobolev classes for the case of special classes of matrices $\{A_{j,k}\}$. Finally, in work \cite{UmarovNew}, which is a continuation of work \cite{UmarovToAppear}, the author obtained a representation of the solution and proved the existence of a weak  solution in the special classes of distributions in the general case, i.e. $\beta$ takes different real values for each equation of the system and the matrix $\{A_{j,k}\}$ is not necessarily triangular. Note that the results
obtained in this paper are new even for time-fractional systems of linear ordinary differential equations.

We will present some comparisons of our results with the already known results of works \cite{Kochuei1} and \cite{UmarovToAppear} after the formulation of the main results (see Remarks \ref{nessesary}, \ref{Umarov}, \ref{Kochubei}, \ref{nonlinearNessesary} and \ref{nonlinear}). But we note right away that in the known works nonlinear systems of fractional partial differential equations were not considered.

In this paper we use the method developed in the works of S. Umarov \cite{UmarovToAppear} and \cite{UmarovNew}, and based on the sequential application of the Fourier and Laplace transforms and the use of some properties of the eigenvalues of Hermitian matrices. Note that when we use this method the differential operators $A_{j,k}(D)$ can have different orders. And the reasoning of the work \cite{Kochuei1} is based on the following two facts: 1) the Green's function of the Cauchy problem for a parabolic system of differential equations was studied in detail in the fundamental monograph by S.D. Eidelman \cite{Eidelman} (see also the monographs of S.G. Krein \cite{Krein}). It should be noted that in the study of the Green's function, an important role was played by the estimate of the norm of the exponential matrix $e^{tA}$ (see Lemma 3.1, \cite{Eidelman}), established in the monograph I.M. Gelfand and G.E. Shilov \cite{GelShil} (see p. 78), which, in turn, is based on the construction  of $e^{tA}$ in the form of some polynomial $R(A)$ (see the monograph by F.R. Gantmacher \cite{Gantmacher}, p. 84). 2) On the other hand, the subordination theorem (see, for example, E. G. Bazhlekova \cite{Bez}) establishes a connection between the Green's function of a parabolic system and a system of fractional differential equations (\ref{system_main}). It should be emphasized that in the monograph \cite{Eidelman}, when establishing the properties of the Green's function, it was required that all operators $A_{j,k}(D)$ have the same order. Therefore, in the paper \cite{Kochuei1}, being based on work \cite{Eidelman}, the main results were obtained when this condition was met.

The paper is organized as follows. Section 2  provides some preliminary facts on the Fourier and Laplace transforms, on the matrix differential operators, and on fractional calculus used in this paper. In section 3 we present results on the existence and uniqueness of a solution to the linear problem (i.e., function $h_j$ does not depend on $u_j$), and in the next section we present the similar results on the nonlinear problem. The article ends with the section Conclusions, which discusses possible generalizations of the results obtained.


\section{Preliminaries and auxiliaries} 
\label{Preliminaries}
In this section we introduce some auxiliary concepts and statements. We recall some properties of the Mittag-Leffler functions, the Laplace transform and estimates of the eigenvalues of Hermitian matrices.

\subsection{The Fourier and Laplace transforms. The Mittag-Leffler functions }

Let us denote by $F[f](\xi)$ (or $\hat{f}(\xi)$) for a given
function $f(x)$ its Fourier transform, and by $F^{-1}f$ the inverse Fourier transform:
\[
F[f](\xi)=\hat{f}(\xi) = \int\limits_{\re^n} f(x) e^{i x \xi} dx, \quad \xi \in \re^n,
\]
and
\[
F^{-1}[\hat{f}](x) =  f(x) = \frac{1}{(2 \pi)^n} \int\limits_{\re^n} \hat{f}(\xi) e^{-i x \xi} d \xi, \quad x \in \re^n.
\]
If, for example, $f\in L_2(\mathbb{R}^n)$, then $F[F^{-1}[f]](x)= F^{-1}[F[f]](x)=f(x)$.

The Laplace transform of a function $\varphi(t)$ of a real variable $t\in \mathbb{R}_+$ is defined by
\begin{equation}\label{laplace}
L[\varphi](s)=\int\limits_0^\infty e^{-st} \varphi (t) dt, \quad s\in \mathbb{C}.
\end{equation}
If this integral  is convergent at the point $s_0\in \mathbb{C}$, then it converges absolutely for $s\in \mathbb{C}$ such that $\Re(s)> \Re(s_0)$. Let us denote  by $\sigma_\varphi$ the infimum of values $\Re (s)$ for which integral (\ref{laplace}) converges.
Then the inverse Laplace transform is given  by the formula
\[
L^{-1}[g](y)=\frac{1}{2\pi i}\int\limits_{\gamma-i\infty}^{\gamma+i\infty} e^{sy} g(s) ds, \quad \gamma >\sigma_\varphi.
\]
Like the Fourier transform, the direct and inverse Laplace transforms are inverses of each other for the “good enough” functions $\varphi$: $L[L^{-1}[\varphi]](t)= L^{-1}[L[\varphi]](t)=\varphi(t)$.

For $0 < \rho < 1$ and an arbitrary complex number $\mu$, let $E_{\rho, \mu}(z)$ denote the Mittag-Leffler function with two parameters of the complex argument $z$ (see, e.g. \cite{KST}, p. 42):
\begin{equation*}
E_{\rho, \mu}(z)= \sum\limits_{k=0}^\infty \frac{z^k}{\Gamma(\rho
k+\mu)}.
\end{equation*}
If the parameter $\mu =1$, then we have the classical
Mittag-Leffler function: $ E_{\rho}(z)= E_{\rho, 1}(z)$.

Recall some properties of the Mittag-Leffler functions.

\begin{lemma}\label{MLestimate} (see \cite{KST}, p. 43). For any $t\geq 0$ one has
\begin{equation*}
|E_{\rho, \mu}(-t)|\leq \frac{C}{1+t},
\end{equation*}
where constant $C$ does not depend on $t$ and $\mu$.
\end{lemma}

\begin{lemma}\label{MLestimate1} (see \cite{KST}, p. 43). For any $t\geq 0$ one has
\begin{equation*}
E_{\rho, \mu}(t) =\frac{1}{\rho}t^{\frac{1}{\rho}(1-\mu)} e^{t^{\frac{1}{\rho}}} +O(t^{-1}).
\end{equation*}
\end{lemma}

Lemma \ref{MLestimate} easily implies the following rougher estimate:
\begin{lemma}\label{MLestimate2} Let $\lambda>0$ and $0<\varepsilon< 1$. Then there is a constant $C>0$, such that
\begin{equation*}
|t^{\rho-1}E_{\rho, \rho}(-\lambda t^\rho)|\leq \frac{C t^{\rho-1}}{1+\lambda t^\rho}\leq C \lambda^{\varepsilon-1} t^{\varepsilon\rho-1}, \quad t>0. 
\end{equation*}
\end{lemma}
\textit{Proof.} Let $\lambda t^{\rho}<1$, then $t<\lambda^{-1/\rho}$ and
\[
t^{\rho-1}=t^{\rho-\varepsilon\rho}t^{\varepsilon\rho-1}<\lambda^{\varepsilon-1}t^{\varepsilon\rho-1}.
\]
If $\lambda t^{\rho}\geq 1$, then $\lambda^{-1}\leq t^{\rho}$ and
\[
\lambda^{-1}t^{-1}=\lambda^{\varepsilon-1}\lambda^{-\varepsilon}t^{-1}\leq \lambda^{\varepsilon-1}t^{\varepsilon\rho-1}.
\]
    \eproof
    
    Let us formulate the following statement with the parameters that we need further, although the statement is known in a more general form.
 \begin{lemma}\label{MLLtransform} (see \cite{KST}, p. 50). Let $t, \mu,$ and $ \lambda>0$. Then
    \[
    L\big[t^{\mu-1}E_{\rho, \mu} (-\lambda t^\rho)\big](s)=\frac{s^{\rho-\mu}}{s^\rho+\lambda}, \,\, \Re(s)>0, \,\, |\lambda s^{-\rho}|<1,
    \]
    and
    \[
    t^{\mu-1}E_{\rho, \mu} (-\lambda t^\rho)= L^{-1}\bigg[\frac{s^{\rho-\mu}}{s^\rho+\lambda}\bigg](t).
    \]
    
\end{lemma}

\subsection{Fractional integrals and derivatives}
Let us recall the definition of fractional derivatives in the sense of Caputo and present their properties that we need. If the function $f(t)$ is locally integrable on $\mathbb{R}_+$, then 
the fractional integral of order $\beta >0$ of function $f$  is defined by
\[
J_t^{\beta} f(t) = \frac{1}{\Gamma(\beta)} \int_0^t (t-\tau)^{\beta-1} f(\tau)d\tau, \, t >0,
\]
where $\Gamma(\beta)$ is Euler's gamma function. 

In this work we will consider only derivatives of order $0 < \beta <1$. The fractional derivative of order $\beta$ of a function $f$ in the sense of Caputo is defined as
\[
D^{\beta}_t f(t)=J^{1-\beta}_t \frac{df}{dt}(t) = \frac{1}{\Gamma(1-\beta)} \int_a^t \frac{f'(\tau)d\tau}{(t-\tau)^{\beta}}, \,\,\, t>0,
\]
provided the integral on the right exists. If we swap the integral with the derivative, then we obtain the definition of the derivative in the Riemann-Liouville sense:

\begin{equation*}
\partial_t^{\beta} f(t) =\frac{d}{dt} \, J_t^{1-\beta} f(t)= \frac{1}{\Gamma(1-\beta)}\frac{d}{dt} \int_a^t \frac{f(\tau)d\tau}{(t-\tau)^{\beta}},
\end{equation*}
provided the expression on the right exists. As an important distinguishing feature of these two derivatives, we note that for equations involving the Caputo derivative, the Cauchy problem with the standard condition $f(0)=f_0$ is well-posed, and for equations involving the Riemann-Liouville derivative, the initial condition should be specified in the form of $J_t^{1-\beta} f(+0)=f_0$.

\begin{lemma}
    
(see \cite{KST}, p. 98)\label{lapl-tr-cap}
The Laplace transform of the
Caputo derivative of a  function $f\in C^1(\mathbb{R}_+)$, 
is
\[
L[D_t^{\beta} f](s)= s^{\beta}L[{f}](s) -
f(0) s^{\beta-1},\,\,\, s>0.
\]
\end{lemma}

We will use this formula in the vector form.   For a vector function $<f_1, \dots, f_m>$ we define 
\[
L[{D}_t^{\beta} <f_1,\dots,f_m>](s)=<L[{D}_t^{\beta}f_1](s), \dots,  L[{D}_t^{\beta} f_m](s)>.
\]
Then
\begin{equation}\label{CD_LT_01} 
L[D_t^{\beta} <f_1,\dots,f_m>](s) = <s^{\beta} L[f_1](s) - f_1(0) s^{\beta-1},\dots,  s^{\beta} L[f_m](s)  - f_m(0) s^{\beta-1}>.
\end{equation}

We also use the following notation for the Riemann-Liouvill fractional derivative of a vector function
$$
\partial_t^\beta <f_1,\dots,f_m>(t)=<\partial_t^{\beta}f_1(t), \dots,  \partial_t^{\beta} f_m(t)>.
$$

\

\subsection{Matrix differential operators with constant symbols}

 Let $A(D)=\sum\limits_{|\alpha|= \ell} a_\alpha D^\alpha$
be a homogeneous differential expression, with constant coefficients and order $l$. Let  $A(\xi)= \sum\limits_{|\alpha|= \ell} a_\alpha \xi^\alpha$
be a symbol of this operator, where $\alpha=(\alpha_1, \alpha_2, ...,
\alpha_n)$ - multi-index and $D=(D_1, D_2, ..., D_n)$,
$D_j=\frac{1}{i}\frac{\partial}{\partial x_j}$, $i=\sqrt{-1}$. 

If $\varphi(x)$ is a sufficiently smooth function, then the action of operator $A(D)$ on this function can be written as:
\begin{equation*}
A(D)\varphi (x)= \frac{1}{(2 \pi)^n} \int_{\re^n} A(\xi) F [\varphi] (\xi)
e^{i x \xi} d \xi\, \quad x \in \re^n.
\end{equation*}
Therefore for such a function one has
\begin{equation}
\label{A(D)}
F[A(D) \varphi(x)] (\xi)=A(\xi) F[\varphi] (\xi).
\end{equation}

Let us recall the definition of an important class of differential operators: operator $A(D)$ is called \textit{elliptic }if inequality $A(\xi)>0$ holds for any $\xi\neq 0$. If $A(D)$ is elliptic, then 
there exist constants $c_1$ and $c_2$ such that
\begin{equation}\label{A}
c_1(1+|\xi|^2)^{\ell} \leq 1+A^{2}(\xi)\leq c_2
(1+|\xi|^2)^{\ell}.
\end{equation}

Let us now consider a matrix of an arbitrary differential expressions $A_{k,j}(D)=\sum_{|\alpha|\leq \ell_{k,j}}a_{\alpha} D^\alpha$ with the order $\ell_{j,k}$ and constant coefficients:
\[
\mathbb{A}(D)= \begin{bmatrix}
	    A_{1,1}(D) & \dots & A_{1, m} (D) \\
	    \dots           & \dots & \dots  \\
	    A_{m, 1}(D) & \dots & A_{m,m} (D)
	    \end{bmatrix}\] 
with the matrix-symbol
\[
\mathcal{A}(\xi) = \begin{bmatrix}
	    A_{1,1}(\xi) & \dots & A_{1, m} (\xi) \\
	    \dots           & \dots & \dots  \\
	    A_{m, 1}(\xi) & \dots & A_{m,m} (\xi)
	    \end{bmatrix},
\]
where $A_{k,j}(\xi)=\sum_{|\alpha|\leq \ell_{k,j}}a_{\alpha} \xi^\alpha$ is a symbol of operator $A_{k,j}(D)$.
Introducing vector-functions $U(t,x)=\langle u_1(t,x),\dots, u_m(t,x) \rangle,$ and 
$H(t,x, U)=\langle h_1(t,x, U),\dots, h_m(t,x, U) \rangle,$ we can represent system (\ref{system_main}) in the vector form:
\begin{equation}
\label{system_1}
D_t^{\beta} {U}(t,x) + \mathbb{A}(D) {U} (t,x) = H(t,x, U),
\end{equation}
where 
$
D_t^{\beta} {U}(t,x)=\langle D_t^{{\beta}} {u}_1(t,x), \dots,  D_t^{{\beta}} {u}_m(t,x) \rangle.
$

Let us list all the conditions for the matrix-symbol $\mathcal{A}(\xi)$ that we will consider fulfilled throughout the entire work. We do not strive for maximum generality and will consider the simplest system of differential equations. Some possible generalizations are given in the Conclusions section. 

\textbf{Conditions $(A)$}. \textit{Assume   that the matrix-symbol $\mathcal{A}(\xi)$ is a Hermitian matrix for all $\xi\in \mathbb{R}^n$ and is positive definite for all $|\xi|\geq R_0$ with a sufficiently large $R_0>0$. We will further assume that all diagonal operators $A_{j,j}(D)$, $j=1, \dots, m$, are homogeneous elliptic differential operators and have the highest order among the operators acting on the function $u_j(t, x)$: }
\begin{equation}\label{maxorder}
    \ell_{j,j}>\ell_{k, j}, \quad \text{for all} \quad k\neq j, \quad j,k=1,\dots, m.
\end{equation}

Let us recall that Hermitianity of $\mathcal{A}(\xi)$ means the fulfillment of equality $A_{k,j}(\xi)=\overline{A_{j,k}(\xi)}$ for all $k,j=1,\dots,m,$ and $\xi\in \mathbb{R}^n$, where $\overline{a}$ is the complex conjugate of $a$ (see e.g., \cite{Horn}, p.7). Positive definiteness of  matrix $\mathcal{A}(\xi)$ means that all eigenvalues of the matrix are positive, or, which is the same, estimate $(\mathcal{A}(\xi)\mu, \mu)>0$,  is valid for all $\mu\in \mathbb{C}^n$ (see e.g., \cite{Bhatia}, p.12).  For the system (\ref{system_main}) a little stricter condition is also known as the uniform strong parabolicity condition (see, e.g., \cite{Handbook}, p. 147, for classical (differential) parabolic systems of equations see the fundamental work by I. Petrowsky \cite{Petrowsky}).

Matrix $\mathcal{A}(\xi)$, like any Hermitian matrix, according to the spectral decomposition theorem, can be represented as a diagonal matrix $\Lambda(\xi)$ translated into another coordinate system (that is, $\mathcal{A}(\xi)=M(\xi)\Lambda(\xi)M^{-1}(\xi) $, where $M(\xi)$ is a unitary matrix, the rows of which are the orthonormal eigenvectors of $\mathcal{A}(\xi)$ that form the basis) (see \cite{Bellman}, Chapter 4). Thus, if $\mathcal{A}(\xi)$ is a positive definite matrix, then all elements of the main diagonal of $\Lambda(\xi)$  (or, in other words, the eigenvalues of $\mathcal{A}(\xi)$) are positive (see, e.g., \cite{Bhatia}, Chapter 1). 

To find the eigenvalues of $\mathcal{A}(\xi)$, we need to solve the following equation for $\lambda$:
\[
|\mathcal{A}(\xi)-\lambda E| =\det \begin{bmatrix}
	    A_{1,1}(\xi)-\lambda & \dots & A_{1, m} (\xi) \\
	    \dots           & \dots & \dots  \\
	    A_{m, 1}(\xi) & \dots & A_{m,m}(\xi)-\lambda 
	    \end{bmatrix} =0,
\]
where $E$ is the identity matrix. This equation, using Vieta's formulas, can be written in the form of algebraic equation
\[
\lambda^m-\lambda^{m-1} (A_{1,1}(\xi)+A_{2,2}(\xi)+\cdots+A_{m,m}(\xi))+\cdots +(-1)^m A_{1,1}(\xi)A_{2,2}(\xi)\cdots A_{m,m}(\xi)=0.
\]
There are $m$ such eigenvalues $\lambda_k(\xi)$ taking into account multiplicity, but there is no exact formula for solutions to this algebraic equation. However using the Gershgorin circle theorem (see, e.g., \cite{Horn}, Theorem 6.1.1), 
  if Conditions $(A)$ are met, we can write an asymptotic estimate for the eigenvalues of the matrix $\mathcal{A}(\xi)$.
Before doing this, note that due to Condition $(A)$ on the the matrix-symbol $\mathcal{A}(\xi)$, one has 
\[
 \lambda_k(\xi)> 0, \quad k=1, \dots, m,\quad  \text{for all} \quad |\xi|>R_0.
   \]

For the convenience of the reader, we formulate and prove the Gershgorin circle theorem  in a form convenient for us.

Let $A$ be a complex $m\times m$ Hermitian matrix with entries $\{a_{j, k}\}$. Let $\varepsilon_j(A)$ be the sum of the absolute values of the non-diagonal entries in the $j$-row
\[
\varepsilon_j(A)=\sum_{k\neq j}|a_{j,k}|. 
\]
Consider the Gershgorin segments
\[
\{t\in \mathbb{R}: \,\, |t-a_{j,j}|\leq \varepsilon_j(A)\}.
\]
\begin{theorem}(The Gershgorin circle theorem). The eigenvalues of the  Hermitian matrix $A$ are in the union of the Gershgorin segments 
\[
G(A)=\bigcup\limits_{j=1}^m \{t\in \mathbb{R}: \,\, |t-a_{j,j}|\leq \varepsilon_j(A)\}
\]
Furthermore, if the union of $k$ of the $m$ segments that comprise $G(A)$ forms a set $G_k(A)$ that
is disjoint from the remaining $m - k$ segmants, then $G_k(A)$ contains exactly $ k$ eigenvalues
of $A$, counted according to their algebraic multiplicities.    
\end{theorem}

\emph{\bf Proof.} Let $\lambda$  be an eigenvalue of $A$ with corresponding eigenvector $\mu= \langle{ \mu_1, \dots, \mu_m\rangle}$. Let $j$ be an index such that $|\mu_j|= \max_{1\leq k\leq m} |\mu_k|$. Since $A\mu=\lambda \mu$, then in particular for the $j$ th component of that equation we get:
\[
\sum_k a_{j,k}\mu_k =\lambda \mu_j.
\]
    Taking $a_{j,j}$ to the other side
 \[
 \sum_{k\neq j} a_{j,k}\mu_k =(\lambda-a_{j,j}) \mu_j.
 \] 
 Therefore, applying the triangle inequality and remembering that $\mu_k\leq \mu_j$, based on how we picked $j$,
 \[
 |\lambda-a_{j,j}|=\bigg|\sum_{k\neq j} a_{j,k}\frac{\mu_k}{\mu_j}\bigg|\leq \sum_{k\neq j} |a_{j,k}|=\varepsilon_j(A),
 \]
  that is, $\lambda$ is in the Gershgorin segmant $\{t\in \mathbb{R}: \,\, |t-a_{j,j}|\leq \varepsilon_j(A)\}$. In particular, $\lambda$ is in the larger set $G(A)$ defined above.

  Now suppose that $k$ of the $m$ segments that comprise $G(A)$ are disjoint from all of the
remaining $m - k$ segments. After a suitable permutation similarity of $A$, we may assume
that $G_k(A) = \sum\limits_{j=1}^k\{t \in \mathbb{R}: |t-a_{j, j}| \leq \varepsilon_j(A)\}$ is disjoint from $G'_k(A) = \sum\limits_{j=k+1}^m\{t \in \mathbb{R}: |t-a_{j, j}| \leq \varepsilon_j(A)\}$. 

 Let $\Lambda$ be the diagonal matrix with entries equal to the diagonal entries of $A$ and let
 \[
 B(t)=(1-t)\Lambda  + tA, \quad 0\leq t\leq 1.
 \]
Observe that $B(0)=\Lambda$, and $B(1)=A$. The diagonal entries of 
$B(t)$ are equal to that of $A$, thus the centers of the Gershgorin segments are the same, however their radii are $t$ times that of $A$: $\varepsilon_j(B(t))=t\varepsilon_j(A)$. Consequently, each of the $m$ Gershgorin  segments of $B(t)$ is contained in the corresponding
Gershgorin  segments  of $A$. In particular
\[
G_k(B(t)) = \sum\limits_{j=1}^k\{t \in \mathbb{R}: |t-a_{j, j}| \leq t \varepsilon_j(A)\}
\]
is contained in $G_k(A)$ and is disjoint from $G'_k(A)$.
Obviously all the eigenvalues $\lambda(t)$
of $B(t)$ are contained in $G(B(t))$, which is contained in $G_k(A) \bigcup G_k'(A)$; $G_k(B(t))$ is disjoint from $G'_k(B(t))$. The segments are closed, so the distance $d$ between two unions $G_k(A)$ and $G'_k(A)$ is positive.

The statement of the theorem is true for $\Lambda=B(0)$ and the distance between unions $G_k(B(t))$ and $G'_k(B(t))$ is a decreasing function of $t$, so it is always at least $d$. Let some eigenvalue $\lambda(0)$ of $\Lambda$ belongs to $ G_k(B(t))$ and since the eigenvalues $\lambda (t)$ of $B(t)$ are a continuous function of $t$ (see e.g., \cite{Voevodin}, Theorem 11.1), then $\lambda(t) \in G_k(B(t))$ for some $t$ and  its distance $d(t)$ from $G'_k(B(t))$ is also continuous. Obviously $d(0)\geq d$, and assume that corresponding eigenvalue  $\lambda (1)$ of $A$ lies in $G'_k(A)$. Then $d(1)=0$, so there exists  $0<t_{0}<1$  such that $ 0<d(t_{0})<d$. But this means $\lambda (t_{0})$ lies outside the Gershgorin segments, which is impossible. Therefore  $\lambda (1)$ lies in the union of the k discs, and the theorem is proven.

   This theorem implies the following statement, which is important for what follows.

   \begin{cor}
       Let the matrix-symbol $\mathcal{A}(\xi)$  satisfy conditions $(A)$. Then the eigenvalues of this matrix can be renumbered  in such a way that the following asymptotic estimates are valid:
       \begin{equation}\label{lambdaEstimate}
  \lambda_j(\xi)=A_{j,j}(\xi)(1+o(1)), \,\,j=1, \dots, m,\,\,\, \text{for all} \,\,\, |\xi|>R_0.
\end{equation}
   \end{cor}
\emph{\bf Proof.}
    It follows from the theorem that at each point $\xi \in \mathbb{R}^n$ with $|\xi|>R_0$, the eigenvalues have estimates 
    \[
    A_{j,j}(\xi) - \ \sum_{k\neq j} |A_{j,k}(\xi)|\leq \lambda_j(\xi)\leq  A_{j,j}(\xi) + \ \sum_{k\neq j} |A_{j,k}(\xi)|
    \]
From here, by virtue of Conditions $(A)$, the statement of the corollary follows.

Next, solving equations $(\mathcal{A}(\xi)-\lambda_j E)\mu$, $j=1,\dots, m,$ with respect to $\mu$, we find  orthonormal eigenvectors $\mu_j(\xi)=\langle{\mu_{j,1}(\xi)\dots, \mu_{j,m}(\xi)\rangle}$, $j=1, \dots, m$, of matrix $\mathcal{A}(\xi)$. Then the matrix $M(\xi)$ will have the form
\[
M(\xi) = \begin{bmatrix}
	    \mu_{1,1}(\xi) & \dots & \mu_{1, m} (\xi) \\
	    \dots           & \dots & \dots  \\
	    \mu_{m, 1}(\xi) & \dots & \mu_{m,m}(\xi) 
	    \end{bmatrix},
\]
with entries $\mu_j(\xi)$: $|\mu_j(\xi)|^2=\sum_{k=1}^m \mu_{j,k}^2(\xi) =1$. 

Let us denote the entries of the inverse matrix $M^{-1}(\xi)$ by $\nu_{j,k}(\xi)$ and note important estimates for further use:
\begin{equation}\label{munu}
 |\mu_{j,k}(\xi)|\leq \mu_0, \quad |\nu_{j,k}(\xi)|\leq \nu_0, \quad \xi\in \mathbb{R}^n, \quad j,k =1, \dots, m.   
\end{equation}

Thus, for matrix $\mathcal{A}(\xi)$, there exists an invertible $(m\times m)$-matrix-function $M(\xi),$ such that
\begin{equation}
\label{matrix_01}
\mathcal{A}(\xi) = M(\xi) \Lambda(\xi) M^{-1}(\xi), \quad \xi\in \mathbb{R}^n,
\end{equation}
with a diagonal matrix 

\[
{\Lambda}(\xi) = \begin{bmatrix}
	    \lambda_1(\xi) & \dots & 0  \\
	    \dots           & \dots & \dots  \\
	    0  & \dots & \lambda_m (\xi)
	    \end{bmatrix}.
\]


\section{Linear system of differential equations}

In this section we investigate the initial-boundary
value problem for the linear (that is  $H(t,x,U)\equiv H(t,x)$) system of differential equation (\ref{system_main}) (see (\ref{system_1})). We emphasize that the reasoning presented here is essentially based on the methods proposed in work S. R. Umarov et al. \cite{UAChen}.

Let $a$ be any positive number and $L_2^a(\mathbb{R}^n)$ stand for the Sobolev space with the scalar product
\[
(\varphi, \psi)_{L_2^a(\mathbb{R}^n)}=\int_{\re^n} F [\varphi] (\xi)\overline{F[\psi]}(\xi) (1+|\xi|^2)^a
 d \xi,\quad \varphi, \,\, \psi \in L_2^a(\mathbb{R}^n)
\]
and the norm
\[
||\varphi||^2_{L_2^a(\mathbb{R}^n)}=\int_{\re^n} |F [\varphi] (\xi)|^2 (1+|\xi|^2)^a
 d \xi, \quad \varphi  \in L_2^a(\mathbb{R}^n).
\]

Further, for any Banach space $B$ we denote by ${\mathbf B}$ the $m$-times topological direct product  
\[
{\mathbf B}  = B \otimes \dots \otimes B,
\] 
of spaces $B$. Elements of ${\mathbf B}$ are vector-functions $\Phi (x) = \langle{ \varphi_1(x), \dots, \varphi_m(x)} \rangle,$ where $\varphi_j(x) \in B, j=1,\dots,m.$ We also introduce the notation 
\[
||\Phi||^2_{\mathbf B}=||\varphi_1||_B^2+\cdots+||\varphi_m||_B^2\,.
\]

As usual, the symbol ${\mathbf C} ([0,T]; {\mathbf B})$ denotes the space of vector-functions continuous on $[0,T]$ with values in ${\mathbf B}$.

\textbf{Linear Problem}. \textit{Find  functions $u_j(t,x)\in L_2^{\ell_{j,j}}(\mathbb{R}^n)$,
$t\in (0, T]$, $j=1, \dots, m$ (note that this inclusion is considered
as a boundary condition at infinity), such that}
\begin{equation}\label{conditionU}
 {U}(t,x)\in {\mathbf C}([0,T]\times \re^n
),  \quad  
 D_t^{\beta} {U}(t,x)\,\,\,\text{and} \,\,\, \mathbb{A}(D) {U} (t,x) \in {\mathbf C}((0,T]\times \re^n),   
\end{equation}
\textit{and satisfying the Cauchy problem}
\begin{align} \label{Cauchy_01_h}
D_t^{\beta} {U}(t,x) &+ \mathbb{A}(D) {U} (t,x)=H(x,t), \quad t>0, \ x \in \re^n,
\\
U(0,x) & =\varPhi(x), \quad x \in \re^n, \label{Cauchy_02_h}
\end{align}
\textit{where $\varPhi(x)$ is a given continuous function.}

Note that the solution to an initial boundary value problem from class (\ref{conditionU}) is usually called \textit{classical solution}.

We draw attention to the fact, that in the statement of the
formulated problem the requirement $u_j(t, x)\in L_2^{\ell_{j,j}}(\mathbb{R}^N)$
is not caused by the merits. However, on the one hand, the
uniqueness of just such a solution is proved quite simply, and on
the other, the solution found by the Fourier method satisfies the
above condition. 

Let us first consider the initial-boundary
value problem for the homogeneous (that is all $h_j\equiv 0$) system of differential equation (\ref{Cauchy_01_h}).
In order to prove the existence and uniqueness of the classical solution (hereinafter simply the solution) of Linear Problem, we proceed as follows: first applying the Fourier and Laplace transforms to the problem, without worrying about the existence of the corresponding integrals, we construct a formal solution. Then we will justify all the reasoning carried out.

So applying Fourier transform to (\ref{Cauchy_01_h}) (with $h_j(t,x)\equiv 0$) and (\ref{Cauchy_02_h}), we obtain a system of fractional order ordinary differential equation (see (\ref{A(D)})):
\begin{equation}\label{euqationFourier}
    D_t^{\beta} F[U](t,\xi) +\mathcal{A}(\xi) F[U](t,\xi)=0, \quad t>0, \ \xi \in \mathbb{R}^n,
\end{equation}
and the initial condition
\begin{equation}\label{conditionFourier}
F[U](0,\xi) = F[\varPhi](\xi), \quad \xi \in \mathbb{R}^n.
\end{equation}
Thus, to determine the unknown function $F[U](t,\cdot)$, we obtained the Cauchy problem (\ref{euqationFourier}), (\ref{conditionFourier}) for a system of ordinary differential equations with constant coefficients (with respect to $t$). This type of Cauchy problem has been studied by many authors (see, for example, \cite{Odabat}, where a representation of the solution is provided). For the convenience of readers, we derive a formula for solving this problem.

Apply the Laplace transform in the vector form \eqref{CD_LT_01}, to get
\begin{align*}
\langle s^{\beta} LF[u_1](s,\xi), \dots, s^{\beta} LF[u_m](s,\xi)  \rangle &+ \mathcal{A}(\xi) LF[U](s,\xi) \\
&= \langle s^{\beta-1} F[\varphi_1](\xi), \dots, s^{\beta-1}F[\varphi_m](\xi)\rangle , \quad s>0, \ \xi \in \mathbb{R}^n.
\end{align*}

Recall that if there is the same number on the main diagonal, and all other elements are equal to zero, then such a matrix is called scalar. Note that the main property of such matrices is that every scalar matrix commutes with any matrix.

Considering this fact and using the representation of matrix $\mathcal{A}(\xi)$ in the form \eqref{matrix_01}, the latter equality can be rewritten as
\[
M (Is^{\beta}+\Lambda(\xi))M^{-1}(\xi) LF[U](s,\xi)= I s^{\beta-1} F[ \varPhi](\xi),
\]
where $I s^{\beta}$, $I s^{\beta-1}$ are scalar matrices with diagonal entries $s^{\beta}, \ s^{\beta-1},$ respectively. From here, transferring the matrices to the right side of the equality, we get (note that it is sufficient to consider case $|\lambda_j(\xi) s^{-\beta}|<1$, $j=1, 2,\dots m,$ see e.g., \cite{KST} p. 50)
\begin{equation}
\label{sol_01}
LF[U](s,\xi) = M(\xi) \mathcal{N}(s,\xi) M^{-1}(\xi) F[\varPhi ](\xi),
\end{equation}
where 
\begin{equation} \label{matrix_03}
{\mathcal{N}}(s,\xi) = \begin{bmatrix}
	   \frac{s^{\beta-1}}{s^{\beta} + \lambda_1(\xi) }& \dots & 0  \\
	    \dots           & \dots & \dots  \\
	    0  & \dots & \frac{s^{\beta-1}}{s^{\beta}+\lambda_m (\xi)}
	    \end{bmatrix}.
\end{equation}
Applying the inverse Laplace transform to equality \eqref{sol_01}, then taking into account Lemma \ref{MLLtransform} and definition \eqref{matrix_03} we obtain the unique solution of the Cauchy problem (\ref{euqationFourier}), (\ref{conditionFourier}):
\[
F[U](t,\xi) = M(\xi) \mathcal{E}_{\beta}(-\Lambda(\xi) t^{\beta}) M^{-1}(\xi)F[\varPhi](\xi), \quad t>0, \ \xi \in \mathbb{R}^n.
\]
Here $\mathcal{E}_{\beta}(-\Lambda(\xi) t^{\beta})$ is the diagonal matrix of the form
\begin{equation} \label{matrix_04}
\mathcal{E}_{\beta}(-\Lambda(\xi) t^{\beta}) = \begin{bmatrix}
	    E_{\beta}(-\lambda_1(\xi)t^{\beta}) & \dots & 0  \\
	    \dots           & \dots & \dots  \\
	    0  & \dots & E_{\beta}(-\lambda_m (\xi)t^{\beta})
	    \end{bmatrix}.
\end{equation}

Thus, the formal  solution of problem \eqref{Cauchy_01_h}-\eqref{Cauchy_02_h} has the representation
\begin{equation}
\label{solution_0001}
U(t,x) = S(t,D) \varPhi(x), \quad t>0, \ x \in \re^n,
\end{equation}
where $S(t,D)$ stands for the  pseudo-differential matrix operator with the matrix-symbol
\[
\mathcal{S}(t,\xi) =  M(\xi) \ \mathcal{E}_{\beta}\Big(- \Lambda(\xi) t^{\beta} \Big)\ M^{-1}(\xi), \quad t>0, \ \xi \in \mathbb{R}^n,
\]
whose entries are
\[
s_{j,k}(t,\xi) = \sum_{q =1}^m \mu_{j,q} (\xi) \nu_{q,k}(\xi) E_{\beta}(-\lambda_{q}(\xi) t^{\beta}), \quad j,k=1,\dots, m,
\]
The explicit component-wise form of the formal solution is
\begin{align}\label{sol_h_0}
u_j(t,x) &= \sum_{k=1}^m s_{j,k}(t,D) \varphi_k(x) \notag 
\\
&= \frac{1}{(2\pi)^n} \sum_{k=1}^m  \sum_{q=1}^m \int\limits_{\re^n} e^{i\xi x}\mu_{j, q}(\xi)E_{\beta} \left( -\lambda_{q}(\xi)t^{\beta} \right) \nu_{q, k}(\xi) F[ \varphi_k] (\xi) d\xi.
\end{align}

Let us denote (see (\ref{lambdaEstimate}))
\begin{equation}\label{minell}
    \tau^* = \ell^{*}- \ell_{*},\quad \ell^{*}=\max_{1\leq k\leq m} \ell_{k,k},\quad  \ell_{*}=\min_{1\leq k\leq m} \ell_{k,k}.
\end{equation}
Note that if all operators $A_{j,j}(D)$ have the same order, then $\tau^*=0$.

Let us also introduce the pseudo-differential operator  $S'(\eta, D)$ with the matrix-symbol 
\[
\mathcal{S}'(t,\xi) =  M(\xi) \ t^{\beta-1}\mathcal{E}_{\beta, \beta}\Big(- \Lambda(\xi) t^{\beta} \Big)\ M^{-1}(\xi), \,\, t>0, \,\, \xi \in \mathbb{R}^n,
\]
and the matrix $\mathcal{E}_{\beta, \beta}$ is defined similarly to (\ref{matrix_04}):
\begin{equation*}
\mathcal{E}_{\beta, \beta}(-\Lambda(\xi) t^{\beta}) = \begin{bmatrix}
	    E_{\beta,\beta}(-\lambda_1(\xi)t^{\beta}) & \dots & 0  \\
	    \dots           & \dots & \dots  \\
	    0  & \dots & E_{\beta, \beta}(-\lambda_m (\xi)t^{\beta})
	    \end{bmatrix}.
\end{equation*}

The following theorem establishes that the formally constructed solution (\ref{sol_h_0}) actually satisfies all conditions (\ref{conditionU}) necessary for the classical solution of problem \eqref{Cauchy_01_h} (with $h_j\equiv 0$) and \eqref{Cauchy_02_h}. A solution to the Linear Problem for an inhomogeneous equation is also constructed.

\begin{theorem}
\label{tfph} Let matrix $\mathcal{A}(\xi)$ satisfy Conditions (A). Let initial and source functions satisfy conditions 
\begin{equation}\label{Con_varphi}
 \varPhi(x) \in {\mathbf L}^{\tau+\tau^*}_2(\mathbb{R}^n),\quad H(t,x) \in {\mathbf C} ([0,T]; {\mathbf L_2}^{\tau+\tau^*}(\mathbb{R}^n)) \quad  \tau > \frac{n}{2}. 
\end{equation}
Then the Linear Problem has a unique solution $U(t,x)=\langle{u_1(t,x), u_2(t,x), \dots, u_m(t,x)}\rangle$ having the representation
\begin{equation}
\label{solution_1}
U(t,x) = S(t,D)\varPhi(x)+\int\limits_0^t S'(\eta, D) H(t-\eta,x) d\eta, \quad 0\leq t\leq T, \ x\in \re^n,
\end{equation}
where $S(t,D)$ and $S'(\eta, D)$ are the pseudo-differential operators defined above. Moreover, each $u_j$  has the property
\begin{equation}\label{bo}
\lim\limits_{|x|\rightarrow\infty} D^\alpha u_j(t, x)=0,\quad
|\alpha|\leq \ell_{j,j}, \quad t>0,
\end{equation}
and for $t>0$ obeys the coercive type estimate
\begin{align}\label{coercive}
||D^\beta_t u_{j}(t,x)||_{C(\mathbb{R}^n)}&+\sum\limits_{k=1}^m||A_{j,k}(D) u_{k}(t,x)||_{C(\mathbb{R}^n)}\notag \\ &\leq C \sum\limits_{k=1}^m \left(t^{-\beta}||\varphi_k||_{L_2^{{\tau+\tau^*}}(\mathbb{R}^n)}+\max_{t\in [0,T]}||h_k(t, \cdot)||_{L_2^{{\tau+\tau^*}}(\mathbb{R}^n)}\right).
\end{align}

\end{theorem}
\begin{remark}
Note that, by virtue of the Sobolev embedding theorem, when conditions (\ref{Con_varphi}) are met, $\varPhi$ and $H$ are continuous functions of their arguments.
\end{remark}

\emph{\bf Proof.} The representation of a formal solution
\[
U(t,x) = S(t,D)\varPhi(x)+\int\limits_0^t S(\eta, D) \partial_t^{1-\beta}H(t-\eta,x) d\eta, \quad 0\leq t\leq T, \ x\in \re^n,
\]
follows directly from \eqref{solution_0001} and from fractional Duhamel's principle \cite{Duhamel_1,Umarov_book_2015}. Note that the integral here is the Laplace convolution. Taking this into account and using the Laplace transform formula for Mittag-Leffler functions (Lemma \ref{lapl-tr-cap}) and the Riemann-Liouville derivative (see e.g., \cite{KST}, p. 284), it is easy to verify that solution $U(t,x)$ can be rewritten in the form (\ref{solution_1}).

It remains to prove that the formal solution (\ref{solution_1}) of the Linear Problem satisfies all its conditions.

Denote the first and second terms on the right of \eqref{solution_1} by $V(t,x)$ and $W(t,x),$ respectively:
\begin{align}\label{V}
V(t,x) &= S(t,D)\varPhi(x), \quad 0\leq t \leq T, \,\, x \in \re^n,           \\
W(t,x) &= \int\limits_0^t S'(\eta, D)  H(t-\eta,x) d\eta, \quad 0\leq t \leq T, \,\, x \in \re^n.\label{W}
\end{align}
Let $v_j(t,x)$ and $w_j(t,x)$ be the entries of the vector functions $V(t,x)$ and $W(t,x)$, respectively. The validity of the statements of the theorem for functions $v_j(t,x)$ and $w_j(t,x)$ is proved in a similar way to work \cite{AZRn}. Let us give a detailed proof for $w_j(t,x)$; for functions $v_j(t,x)$ is proved using the same reasoning.

Functions $w_j(t,x)$ can be written in the following explicit component-wise form
\begin{equation}\label{sol_fi_0}
w_j(t,x)=  \frac{1}{(2\pi)^n} \sum_{k=1}^m  \sum_{q=1}^m \int\limits_0^t \int\limits_{\mathbb{R}^n} e^{i\xi x}\mu_{j, q}(\xi)\eta^{\beta-1} E_{\beta, \beta} \left( -\lambda_{q}(\xi)\eta^{\beta} \right) \nu_{q, k}(\xi) F[ h_k] (t-\eta, \xi) d\xi d\eta. 
\end{equation}
To investigate functions $w_j(t,x)$, it is sufficient to consider the following integrals (see Lemma \ref{MLestimate1})
\[
W_{j, R}(t,x, k, q)=\int\limits_0^t \int\limits_{R_0<|\xi|<R} e^{i\xi x}a_j(\xi)\eta^{\beta-1}E_{\beta, \beta} \left( -\lambda_q(\xi)\eta^{\beta} \right) \hat{h}_k(t-\eta, \xi) d\xi d\eta,
\]
where $|a_j(\xi)|=|\mu_{j,q}(\xi) \nu_{q,k}(\xi)|\leq C$ (see (\ref{munu})). Here the number $R_0$ is taken from estimate (\ref{lambdaEstimate}) and therefore eigenvalues $\lambda_q(\xi)$ satisfy this estimate.

Suppose that the functions  $\varphi_k(x)$ and $h_k(x,t)$  satisfy all the conditions
of Theorem \ref{tfph} and let $|\alpha|\leq \ell_{j,j}$. Then one has
$$
||D^\alpha W_{j, R}(t,x, k, q)||_{C(\mathbb{R}^n)}=
$$
$$
\bigg|\bigg|\int\limits_0^t\int\limits_{R_0<|\xi|<R} e^{i\xi x}a_j(\xi) \xi^\alpha\eta^{\beta-1}E_{\beta, \beta} \left( -\lambda_q(\xi)\eta^{\beta} \right) \hat{h}_k(t-\eta, \xi) d\eta d\xi \bigg|\bigg|_{C(\mathbb{R}^n)}\leq
$$
$$
\int\limits_0^t \int\limits_{R_0<|\xi|<R} \bigg|a_j(\xi) \xi^\alpha \eta^{\beta-1}E_{\beta, \beta} \left( -\lambda_q(\xi)\eta^{\beta} \right) \hat{h}_k(t-\eta, \xi)\bigg|  d\xi d\eta.
$$
We have $|a_j(\xi)|\leq C$ and therefore Lemma \ref{MLestimate2} and estimates (\ref{A}), (\ref{lambdaEstimate}) imply (see definition (\ref{minell}))
\[
\big|a_j(\xi)\xi^\alpha \eta^{\beta-1} E_{\beta, \beta}(-\lambda_q(\xi)\eta^\beta)\big|\leq C\eta^{\varepsilon\beta-1}\lambda_q^{\varepsilon-1}(\xi) |\xi|^{|\alpha|}\leq C\eta^{\varepsilon\beta-1} |\xi|^{\ell_{j,j}-\ell_{*}(1-\varepsilon)}.
\]
Now set $\tau>n/2$ and choose $\varepsilon>0$ such that $\tau-\varepsilon \ell_* >n/2$. Then the Holder
inequality implies
$$
||D^\alpha W_{j, R}(t,x, k, q)||_{C(\mathbb{R}^n)}\leq C\int\limits_0^t\eta^{\varepsilon\beta-1}d\eta \int\limits_{R_0<|\xi|<R}|\xi|^{\tau+\ell_{j,j}-\ell_{*}}|\hat{h}_k(t-\eta, \xi)| |\xi|^{-\tau+\varepsilon\ell_*}d\xi\leq
$$
\[
C_{\varepsilon, \beta, \tau}\max_t\bigg|\bigg|(1+|\xi|^2)^{\frac{1}{2}\big[\tau+\ell_{j,j}-\ell_{*}\big]}\hat{h}_k(t, \xi)\bigg|\bigg|_{L_2(\mathbb{R}^n)}\leq
C_{\varepsilon, \beta, \tau}\max_{t\in [0,T]}||h_k(t, \cdot)||_{L_2^{{\tau+\ell_{j,j}-\ell_{*}}}(\mathbb{R}^n)}.
\]
Thus one has
\begin{equation}\label{estimateW}
||D^\alpha w_{j}(t,x)||_{C(\mathbb{R}^n)}\leq C \sum\limits_{k=1}^m \max_{t\in [0,T]}||h_k(t, \cdot)||_{L_2^{{\tau+\tau^*}}(\mathbb{R}^n)},
\end{equation}
and therefore $\mathbb{A}(D) {W} (t,x) \in {\mathbf C}([0,T]\times
\re^n)$, and in particular, ${W}(t,x)\in {\mathbf C}([0,T]\times
\re^n)$.

The fact that function $V(t,x)$ has similar properties is proven in the same way. So one has
\[
||V_{j, R}(t,x, k, q)||_{C(\mathbb{R}^n)}\leq C ||\varphi_k||_{L_2^{{\tau}}(\mathbb{R}^n)},
\]
and for $|\alpha|\leq \ell_{j,j}$
\begin{equation}\label{estimateV}
||D^\alpha V_{j, R}(t,x, k, q)||_{C(\mathbb{R}^n)}\leq C t^{-\beta}||\varphi_k||_{L_2^{{\tau+\tau^*}}(\mathbb{R}^n)},
\end{equation}
where
\[
V_{j, R}(t,x, k, q)=\int\limits_{R_0<|\xi|<R} e^{i\xi x}a_j(\xi)E_{\beta} \left( -\lambda_{q}(\xi)t^{\beta} \right) \hat{\varphi_k} (\xi) d\xi.
\]
Therefore $\mathbb{A}(D) {V} (t,x) \in {\mathbf C}([0,T]\times
\re^n)$, and  ${V}(t,x)\in {\mathbf C}([0,T]\times
\re^n)$.

Hence, in view of the equality $U=V+W$ we have,
$$ {U}(t,x)\in {\mathbf C}([0,T]\times
\re^n),  \quad  
 \mathbb{A}(D) {U} (t,x) \in {\mathbf C}((0,T]\times
\re^n).  
$$
Since by virtue of equation (\ref{Cauchy_01_h}) we have $ D_t^{\beta} {U}(t,x) = - \mathbb{A}(D) {U} (t,x) + H(t,x)$, then
$D_t^{\beta} {U}(t,x)\in {\mathbf C}((0,T]\times
\re^n)$.

Similarly the following estimates for $|\alpha|\leq \ell_{j,j}$ are proved:
\[
||D^\alpha W_{j, R}(t,x, k, q)||_{L_2(\mathbb{R}^n)}\leq C ||h_k(t,\cdot)||_{L_2^{{\tau^*}}(\mathbb{R}^n)}\,,
\]
\[
||D^\alpha V_{j, R}(t,x, k, q)||_{L_2(\mathbb{R}^n)}\leq C t^{-\beta}||\varphi_k||_{L_2^{{\tau^*}}(\mathbb{R}^n)}, \quad |\alpha|>0,
\]
and 
\begin{equation}\label{estimateV_0}
  ||V_{j, R}(t,x, k, q)||_{L_2(\mathbb{R}^n)}\leq C ||\varphi_k||_{L_2^{{\tau^*}}(\mathbb{R}^n)}\,.  
\end{equation}
It follows that $u_{j,j}(t,x)\in L_2^{\ell_{j,j}}(\mathbb{R}^n)$, $t\in (0, T]$, $j=1, \dots, m$. 

The conditions $u_{j,j}(t,x)\in L_2^{\ell_{j,j}}(\mathbb{R}^n)$, $t\in (0, T]$, $j=1, \dots, m$, of Linesr Problem imply the uniqueness of the solution. Indeed, consider the homogeneous Linear Problem (i.e. $\varPhi(x) \equiv 0$ and $H(t, x)\equiv 0$). For such functions $u_{j,j}(t,x)$, equation (\ref{euqationFourier}) and the homogeneous condition (\ref{conditionFourier}) are satisfied for almost all $\xi\in \mathbb{R}^n$ and all $t\geq 0$. Consequently, $F[U](t,\xi)=0$ almost everywhere $\xi\in \mathbb{R}^n$ and $t>0$. In view of Parseval's equality $U(t,x)=0$ for almost all $x\in\mathbb{R}^n$ and $t\geq 0$. Due to the continuity of this function, we obtain the uniqueness of the solution.

Let us show the property (\ref{bo}) of the solution. To do this, note first that
the inclusion $\varphi_j\in L_2^{\tau+\tau^*}(\mathbb{R}^n)$ implies $(1+|\xi|^2)^{\tau^*/2}\hat{\varphi}_j(\xi)\in L_1(\mathbb{R}^n)$, $j=1, \dots, m$. Indeed, application of the
Holder inequality gives
\[
\int\limits_{\mathbb{R}^n}(1+|\xi|^2)^{\tau^*/2}|\hat{\varphi}_j(\xi)|d\xi =\int\limits_{\mathbb{R}^n}(1+|\xi|^2)^{(\tau+\tau^*)/2}(1+|\xi|^2)^{-\tau/2}|\hat{\varphi}_j(\xi)|d\xi \leq C_\tau ||\varphi_j||_{L_2^{\tau+\tau^*}(\mathbb{R}^n)}.
\]
Therefore, by estimate (\ref{MLestimate}), one has
\[
\big|\xi^{\alpha}a_j(\xi)  E_{\beta}(-\lambda_q(\xi)t^\beta)\hat{\varphi}_k(\xi)\big|\leq Ct^{-\beta}(1+|\xi|^2)^{\tau^*/2}|\hat{\varphi}_k(\xi)|\in L_1(\mathbb{R}^n)
\]
for all $|\alpha|\leq \ell_{j,j}$ and $t>0$. Hence $\lim\limits_{|x|\rightarrow\infty} D^\alpha V_{j, R}(t,x, k, q)=0,\quad
|\alpha|\leq \ell_{j,j}, \quad t>0$. Functions $W_{j, R}(t,x, k, q)$ have a similar property. Thus, property (\ref{bo}) of the solution is proven.

The estimate (\ref{coercive}) is a consequence of equation (\ref{Cauchy_01_h}) and estimates (\ref{estimateW}) and (\ref{estimateV}).

Theorem \ref{tfph} is completely proven.
\eproof

\begin{remark}\label{nessesary}
If all operators $A_{j,j}(D)$  have the same order, then $\tau^*=0$. Therefore, the sufficient condition for the initial functions and for the right-hand sides of the equations has the form $\varphi_j, \, h_j(t,x) \in L^{\tau}_2(\mathbb{R}^n),$ $\tau>\frac{n}{2}.$ Each such function, by virtue of the Sobolev embedding theorem, is continuous, which is required in the conditions of the problem. Therefore, condition $\tau>\frac{n}{2}$ cannot be weakened, i.e. even if $\tau =\frac{n}{2}$ then in the class $L^{\tau}_2(\mathbb{R}^n)$ there are also unbounded functions.
\end{remark}

\begin{remark}\label{Umarov}
 In work \cite{UmarovToAppear}, when studying the existence of a strong solution from Sobolev classes, it is assumed that the eigenvalues of matrix $\mathcal{A}(\xi)$  satisfy estimates of type (\ref{lambdaEstimate}). The present work shows for which matrices these estimates hold and in Theorem \ref{tfph} sufficient conditions are found (in some cases they are necessary) on the initial function and on the right-hand side of the equation, which ensure the existence of a classical solution.   
\end{remark}

Taking into account formulas (\ref{V}) and (\ref{W}), we will use the definition of the Fourier transform in the integrals (\ref{sol_h_0}) and (\ref{sol_fi_0}), and then change the order of integration in these integrals (note that for the functions under consideration this is possible). Then we will have
\[
v_j(t,x)=\sum\limits_{k=1}^m \int\limits_{\mathbb{R}^n} Z_{\{j,k\}, \beta}(t, x-y) \varphi_k(y) dy,
\]
and
\[
w_j(t,x)=\sum\limits_{k=1}^m \int\limits_0^t\int\limits_{\mathbb{R}^n} Y_{\{j,k\}, \beta}(t-\eta, x-y) h_k(\eta, y) dyd\eta,
\]
where
\[
Z_{\{j,k\}, \beta}(t, x)=\sum\limits_{q=1}^m\int\limits_{\mathbb{R}^n} e^{ix\xi}\mu_{j, q}(\xi) E_{\beta} \left( -\lambda_{q}(\xi)t^{\beta} \right) \nu_{q, k}(\xi)d\xi,
\]
and
\[
Y_{\{j,k\}, \beta}(t, x)=\sum\limits_{q=1}^m\int\limits_{\mathbb{R}^n}e^{ix\xi}\mu_{j, q}(\xi)t^{\beta-1} E_{\beta, \beta} \left( -\lambda_{q}(\xi)t^{\beta} \right) \nu_{q, k}(\xi)d\xi.
\]
Thus, if we introduce matrices $Z_\beta(t,x)=\{Z_{\{j,k\}, \beta}(t, x)\}$ and $Y_\beta(t,x)=\{Y_{\{j,k\}, \beta}(t, x)\}$, then for the solution of Linear Problem we obtain the representation
\begin{equation}\label{solution_LP}
    U(t,x)= \int\limits_{\mathbb{R}^n} Z_\beta(t,x-y)\varPhi(y)dy+\int\limits_0^t\int\limits_{\mathbb{R}^n}Y_\beta(t-\eta,x-y) H(\eta, y) dy d\eta.
\end{equation}
\begin{remark}\label{Kochubei}
    In the case when all operators $A_{j,k}(D)$ have the same order $2b$, the representation (\ref{solution_LP}) for the solution of Linear Problem is established in work \cite{Kochuei1} (see also \cite{Kochubei3} and \cite{Kochubei4}). In this work, the matrix functions $Z_\beta(t,x)$ and $Y_\beta(t,x)$ are defined using Wright functions. When constructing Green's functions, the author used the well-known subordination method (see, for example, \cite{Baz}). The main goal in these works is to obtain estimates for functions $Z_\beta(t,x)$ and $Y_\beta(t,x)$ and prove the uniqueness of the solution to the Cauchy problem. Therefore, in these works the question of finding conditions on the initial function and on the right-hand side of the equation that ensure the existence of a solution to Linear Problem was not explored. Let us note once again that in \cite{Kochubei4} the author also considered the case when the differential expressions $A_{j,k}(D)$ depend on the spatial variables $x$.
\end{remark}
\textbf{Example}.
To illustrate the theorem proved above consider  the following initial-boundary value problem (slightly modified problem from the article \cite{UAChen}): Find the functions $u_1(t,x), u_2(t,x)\in L_2^2(\mathbb{R}^n)$ satisfying the equations and the initial conditions
\begin{align} \label{ex_01}
D_t^{\beta} u_1(t,x) & - \Delta_x  u_1(t,x)+ i\sum\limits_{j=1}^n a_j \frac{\partial u_2}{\partial x_j}  (t,x)=0, \quad t>0, \ x\in \mathbb{R}^n,
\\\label{ex_02}
D_t^{\beta} u_2(t,x) & + i\sum\limits_{j=1}^n a_j \frac{\partial u_1}{\partial x_j}  -\Delta_x  u_2(t,x)=0, \quad t>0, \ x\in \mathbb{R}^n
\\ \label{ex_03}
u_1(0,x) & =\varphi_1(x), \quad u_2(0,x)=\varphi_2(x),  \quad x\in \mathbb{R}^n,
\end{align}
where $\Delta_x$ is the Laplace operator, $a_j$ are real numbers, $\beta_j\in (0,1]$ and $\varphi_j$ are given functions that satisfy the conditions of Linear Problem.

It is not hard to see that the matrix-symbol $\mathcal{A}(\xi)$ of the operator corresponding to equations \eqref{ex_01}-\eqref{ex_02} is Hermitian and has the representation

\begin{equation*} \label{matrix_10}
 \begin{bmatrix}
	    |\xi|^2 & \sum\limits_{j=1}^n a_j\xi_j \\
	    \sum\limits_{j=1}^n a_j\xi_j &  |\xi|^2
	    \end{bmatrix} 
	    =
	    \begin{bmatrix}
	    \sqrt{2}/2 & \sqrt{2}/2 \\
	    -\sqrt{2}/2 &  \sqrt{2}/2
	    \end{bmatrix}
	    \begin{bmatrix}
	    |\xi|^2+\sum\limits_{j=1}^n a_j\xi_j & 0 \\
	    0 &  |\xi|^2-\sum\limits_{j=1}^n a_j\xi_j
	    \end{bmatrix}	
	        \begin{bmatrix}
	    \sqrt{2}/2 & -\sqrt{2}/2 \\
	    \sqrt{2}/2 &  \sqrt{2}/2
	    \end{bmatrix}	.   
	    \end{equation*}
    It follows that the orthonormal eigenvectors of matrix $\mathcal{A}(\xi)$ have the form $\mu_1(\xi)=(\mu_{1,1}(\xi), \mu_{1,2}(\xi))=( \sqrt{2}/2,  -\sqrt{2}/2) $, and $\mu_2(\xi)=(\mu_{2,1}(\xi), \mu_{2,2}(\xi))=( \sqrt{2}/2,  \sqrt{2}/2) $ and the corresponding eigenvalues can be represented in the form  $\lambda_1(\xi)= |\xi|^2+\sum\limits_{j=1}^n a_j\xi_j$ and $\lambda_2(\xi)=|\xi|^2-\sum\limits_{j=1}^n a_j\xi_j.$ Note in this example the eigenvectors does not depend on $\xi$. For the elements $\{\nu_{j,k}\}$ of the inverse matrix $M^{-1}$we have: $\nu_1(\xi)=(\nu_{1,1}(\xi), \nu_{1,2}(\xi))=( \sqrt{2}/2,  \sqrt{2}/2) $, and $\nu_2(\xi)=(\nu_{2,1}(\xi), \nu_{2,2}(\xi))=( -\sqrt{2}/2, \sqrt{2}/2) $.
    
    The symbol of the solution operator $S(t,D)$ is the matrx $\mathcal{S}(t,\xi)=\{s_{j,k}(t,\xi) \}, \ j,k=1,2,$ with entries
\[
s_{1,1} (t,\xi)=s_{2,2} (t,\xi) = \frac{1}{2} E_{\beta}((|\xi|^2+\sum\limits_{j=1}^n a_j\xi_j)t^{\beta})+\frac{1}{2}E_{\beta}((|\xi|^2-\sum\limits_{j=1}^n a_j\xi_j)t^{\beta}),
\]
\[
s_{1,2}(t,\xi) =s_{2,1}(t,\xi) = \frac{1}{2} E_{\beta}((|\xi|^2+\sum\limits_{j=1}^n a_j\xi_j)t^{\beta})-\frac{1}{2}E_{\beta}((|\xi|^2-\sum\limits_{j=1}^n a_j\xi_j)t^{\beta}).
    \]
Therefore, the formal solution   $U(t,x)= \langle u_1(t,x), u_2(t,x) \rangle$ to initial-boundary value problem \eqref{ex_01}-\eqref{ex_03} has the representation:
\begin{align*} \label{sol_11}
u_{1} (t,x)&= \left[\frac{1}{2} E_{\beta}((\Delta_x+i\sum\limits_{j=1}^n a_jD_j)t^{\beta})+\frac{1}{2}E_{\beta}((\Delta_x-i\sum\limits_{j=1}^n a_jD_j)t^{\beta}) \right]\varphi_1(x) 
\\
&+ 
 \left[\frac{1}{2} E_{\beta}((\Delta_x+i\sum\limits_{j=1}^n a_jD_j)t^{\beta})-\frac{1}{2}E_{\beta}((\Delta_x-i\sum\limits_{j=1}^n a_jD_j)t^{\beta}) \right]\varphi_2(x);
\\
u_{2} (t,x)&= \left[\frac{1}{2} E_{\beta}((\Delta_x+i\sum\limits_{j=1}^n a_jD_j)t^{\beta})-\frac{1}{2}E_{\beta}((\Delta_x-i\sum\limits_{j=1}^n a_jD_j)t^{\beta}) \right]\varphi_1(x) 
\\
&+ 
 \left[\frac{1}{2} E_{\beta}((\Delta_x+i\sum\limits_{j=1}^n a_jD_j)t^{\beta})+\frac{1}{2}E_{\beta}((\Delta_x-i\sum\limits_{j=1}^n a_jD_j)t^{\beta}) \right]\varphi_2(x).
	    \end{align*}
	    
	    In the explicit form this solution has the form
\begin{align*} 
u_{1} (t,x)&= \frac{1}{(2\pi)^n} \int\limits_{\mathbb{R}^n} \left[\frac{1}{2} E_{\beta}((|\xi|^2+\sum\limits_{j=1}^n a_j\xi_j)t^{\beta})+\frac{1}{2}E_{\beta}((|\xi|^2-\sum\limits_{j=1}^n a_j\xi_j)t^{\beta}) \right] F[{\varphi}_1](\xi) d \xi 
\\
&+ 
\frac{1}{(2\pi)^n} \int\limits_{\mathbb{R}^n}\left[\frac{1}{2} E_{\beta}((|\xi|^2+\sum\limits_{j=1}^n a_j\xi_j)t^{\beta})-\frac{1}{2}E_{\beta}((|\xi|^2-\sum\limits_{j=1}^n a_j\xi_j)t^{\beta}) \right] F[{\varphi}_2](\xi) d\xi;
\\
u_{2} (t,x)&= \frac{1}{(2\pi)^n} \int\limits_{\mathbb{R}^n} \left[\frac{1}{2} E_{\beta}((|\xi|^2+\sum\limits_{j=1}^n a_j\xi_j)t^{\beta})-\frac{1}{2}E_{\beta}((|\xi|^2-\sum\limits_{j=1}^n a_j\xi_j)t^{\beta}) \right] F[{\varphi}_1](\xi) d \xi 
\\
&+ 
\frac{1}{(2\pi)^n} \int\limits_{\mathbb{R}^n} \left[\frac{1}{2} E_{\beta}((|\xi|^2+\sum\limits_{j=1}^n a_j\xi_j)t^{\beta})+\frac{1}{2}E_{\beta}(|\xi|^2-\sum\limits_{j=1}^n a_j\xi_j)t^{\beta}) \right] F[{\varphi}_2](\xi) d \xi.
	    \end{align*}

It is not hard to see that $\lambda_k(\xi) \geq 0, \ k=1,2,$ for all $\xi$ satisfying the inequality $|\xi| \ge \max_j |a_j|.$ Moreover, if $|\xi|$ is taken even larger, say $|\xi|\geq 2|a_j|$, then $\lambda_k(\xi) \geq 1/2 |\xi|^2, \ k=1,2,$ (see (\ref{lambdaEstimate})). Note that both operators $A_{1,1}(D)$ and $A_{2,2}(D)$ are the Laplace operator, i.e. have second order. Therefore $\ell^{*}=2$  and $\tau^*=0$ (see (\ref{minell})).

Applying Theorem \ref{tfph} we see that
$U(t,x)$  satisfies all the conditions (\ref{conditionU}) if $\varphi_k \in L_2^\tau(\mathbb{R}^n), \, \tau>\frac{n}{2}, \,\, k=1,2.$ By virtue of the Sobolev embedding theorem, all functions from this class are continuous. As noted above, this condition is unimprovable in the sense that if $\tau=\frac{n}{2}$, then in the class $L_2^{\frac{n}{2}}$ there may exist unbounded functions $\varphi$ and the initial conditions will not take place in the classical sense.

\section{Nonlinear system of differential equations}

Now let us consider the initial-boundary
value problem for the nonlinear system of differential equation (\ref{system_main}) (see (\ref{system_1})). We will require the nonlinear function $H(t,x,U)$ to be continuous on $t\in [0.T]$ and $x\in \mathbb{R}^n$ and to satisfy the standard Lipschitz condition for such equations: for any $Y=\langle{y_1, \dots, y_m\rangle}$ and $Z=\langle{z_1, \dots, z_m\rangle}$ we assume, that
\begin{equation}\label{H_L}
    |h_j(t,x, Y)- h_j(t,x, Z)|\leq L_0 \sum\limits_{k=1}^m |y_k-z_k|,\,\, t\in [0, T],\,\, x\in \mathbb{R}^n, \,\, j=1, \cdots, m,
    \end{equation}
where $L$ is some positive constant independent of $t,x, Y, Z$.

It is also important for what follows that the condition (see condition (\ref{Con_varphi}) of Theorem \ref{tfph})
\begin{equation}\label{H0}
    H(t, x, 0)\in  {\mathbf C} ([0,T]; {\mathbf L_2}^{\tau+\tau^*}(\mathbb{R}^n)) ,\quad  \tau > \frac{n}{2},
\end{equation}
is satisfied. If we now denote $_0 H (t, x, U)= H(t,x, U)- H(t, x, 0)$ with the entries $_0 h_j (t,x, U)$, then $_0 H(t,x, 0)=0,$ $ H (t, x, U) =\,\,\, _0 H(t,x, U)+ H(t, x, 0)$ and Lipschitz’s condition (\ref{H_L}) implies
\begin{equation}\label{0_H_L}
    |_0 h_j(t,x, Y)|\leq L_0 \sum\limits_{k=1}^m |y_k|,\,\,  t\in [0, T],\,\, x\in \mathbb{R}^n, \,\, j=1, \cdots, m.
    \end{equation}

\textbf{Nonlinear Problem}. \textit{Find  functions $u_j(t,x)\in L_2^{\ell_{j,j}}(\mathbb{R}^n)$,
$t\in [0, T]$, $j=1, \dots, m$ (note that this inclusion is considered
as a boundary condition at infinity), such that}
\begin{equation}\label{NconditionU}
{U}(t,x)\in {\mathbf C}([0,T]\times \re^n
),  \quad  
 D_t^{\beta} {U}(t,x)\,\,\,\text{and} \,\,\, \mathbb{A}(D) {U} (t,x) \in {\mathbf C}((0,T]\times \re^n),   
\end{equation}
\textit{
and satisfying the Cauchy problem}
\[
D_t^{\beta} {U}(t,x) + \mathbb{A}(D) {U} (t,x)=H(t, x, U), \quad t>0, \ x \in \re^n,
\]
\[
U(0,x)  =\varPhi(x), \quad x \in \re^n, \label{NCauchy_02_h}
\]
\textit{where $\varPhi(x)$ is a given continuous function and the source function $H(t,x, U)$ satisfies conditions }(\ref{H_L}) and (\ref{H0}). 

Note that the solution to an initial boundary value problem from class (\ref{NconditionU}) is usually called \textit{classical solution}.

The following statement is true.

\begin{lemma}\label{LestimateH} Let the vector-functions $Y(t,x), \, Z(t,x) \in {\mathbf C} ([0,T]; {\mathbf L_2}(\mathbb{R}^n))$. Then for any $\eta_1, \, \eta_2 \in [0,T]$ the estimate
\[
    \bigg|\bigg|\int\limits_{\eta_1}^{\eta_2} S'(\eta_2-\eta, D)  \big[H(\eta,x, Y)-H(\eta,x, Z)\big] d\eta\bigg|\bigg|_{{\mathbf L}_2(\mathbb{R}^n)}\leq 
    \]
    \begin{equation}\label{EestimateH}
    c_1 (\eta_2-\eta_1)^\beta \sum\limits_{k=1}^m \max_{t\in [0,T]} ||y_k(t,\cdot)-z_k(t,\cdot)||_{L_2(\mathbb{R}^n)}, 
\end{equation}
is valid. Where the constant $c_1$ depends only on numbers ${L,\beta, T, R_0}$.
\end{lemma}

\emph{\bf Proof.} Again, as in the proof of Theorem \ref{tfph}, to prove the lemma, it is sufficient to consider the following integrals
\[
R_{j}(\eta_1, \eta_2,x, k, q)=
\]
\[
\int\limits_{\eta_1}^{\eta_2} \int\limits_{\mathbb{R}^n} e^{i\xi x}a_j(\xi)(\eta_2-\eta)^{\beta-1}E_{\beta, \beta} \left( -\lambda_q(\xi)(\eta_2-\eta)^{\beta} \right) [\hat{h}_k(\eta, \xi, Y) -\hat{h}_k(\eta, \xi, Z)] d\xi d\eta.
\]
Apply Lemmas \ref{MLestimate} and \ref{MLestimate1} to get (pay attention to the definition of the number $R_0$ in (\ref{lambdaEstimate}))
\begin{equation}\label{estimateKernel}
    |a_j(x) E_{\beta, \beta} \big( -\lambda_q(x)(\eta_2-\eta)^{\beta} \big)|\leq C_{T, R_0}, \,\, x\in \mathbb{R}^n.
\end{equation}
Hence, by virtue of Parseval’s equality, we have 
\[
||R_{j}(\eta_1, \eta_2,x, k, q)||_{L_2(\mathbb{R}^n)}\leq
C_{T, R_0}\int\limits_{\eta_1}^{\eta_2} (\eta_2-\eta)^{\beta-1}||h_k(\eta, \cdot, Y) -h_k(\eta, \cdot,  Z) ||_{L_2(\mathbb{R}^n)} d\eta.
\]
 Finally, using the Lipschitz condition (\ref{H_L}) we obtain the required estimate (\ref{EestimateH}).

\begin{lemma}\label{WestimateH} Let a vector-functions $Y(t,x) \in {\mathbf C} ([0,T]; {\mathbf L_2}(\mathbb{R}^n))$. Then the estimate
\begin{equation}\label{EestimateW}    \bigg|\bigg|\int\limits_{0}^{t} S'(t-\eta, D)  H(\eta,x, Y) d\eta\bigg|\bigg|_{{\mathbf L}_2(\mathbb{R}^n)}\leq 
     c_2  \sum\limits_{k=1}^m \max_{t\in [0,T]}\left( ||y_k(t,\cdot)||_{L_2(\mathbb{R}^n)}+||h_k(t,\cdot, 0)||_{L_2(\mathbb{R}^n)}\right), 
\end{equation}
is valid. Where the constant $c_2$ depends only on numbers ${L,\beta, T, R_0}$.
\end{lemma}
\emph{\bf Proof.}
Set    
\[
W_{j}(t,x, k, q)=
\]
\[
\int\limits_{0}^{t} \int\limits_{\mathbb{R}^n} e^{i\xi x}a_j(\xi)(t-\eta)^{\beta-1}E_{\beta, \beta} \left( -\lambda_q(\xi)(t-\eta)^{\beta} \right) \hat{h}_k(\eta, \xi, Y) d\xi d\eta.
\]
Then the estimate of the kernel (\ref{estimateKernel}) implies
  \[
||W_{j}(t,x, k, q)||_{L_2(\mathbb{R}^n)}\leq
C_{T, R_0}\int\limits_{0}^{t} (t-\eta)^{\beta-1}||_0 h_k(\eta, \cdot, Y) + h_k(\eta, \cdot,  0) ||_{L_2(\mathbb{R}^n)} d\eta.
\] 
Apply conditions (\ref{H0}) and (\ref{0_H_L}) to get the required estimate (\ref{EestimateW}).

Let a vector-function $F(t,x)=\langle{f_1(t,x), \cdots, f_m(t,x)\rangle}\in {\mathbf C} ([0,T]; {\mathbf L_2}(\mathbb{R}^n))$. Consider the Volterra integral equation
\begin{equation}\label{Volterra}
 U(t,x) = F(t,x)+\int\limits_0^t S'(t-\eta, D) H(\eta,x, U) d\eta, \quad 0\leq t\leq T, \ x\in \re^n.
\end{equation}

\begin{lemma}\label{VolterraSolution}
   There exists a unique solution $U(t,x) \in  {\mathbf C} ([0,T]; {\mathbf L_2}(\mathbb{R}^n))$ to the integral equation (\ref{Volterra}).
\end{lemma}

\emph{\bf Proof.}
    Equation (\ref{Volterra}) is similar to the equations considered in the book \cite{KST}, p. 199, Eq. (3.5.4),
and it is solved in essentially the same way. 

Equation (\ref{Volterra}) makes sense in any interval $[0, t_1] \subset  [0, T ], 0 < t_1 < T$. Choose $t_1$ such that
\begin{equation}\label{c1}
    c_1 t_1^\beta<1,
\end{equation}
where constant $c_1$ is from estimate (\ref{EestimateH}).
Let us prove the existence of a unique solution $U(t,x) \in {\mathbf C} ([0,t_1]; {\mathbf L_2}(\mathbb{R}^n))$  to the equation (\ref{Volterra})
on the interval $[0, t_1]$. For this we use the Banach fixed point theorem for the space ${\mathbf C} ([0,t_1]; {\mathbf L_2}(\mathbb{R}^n))$ (see, e.g., \cite{KST}, theorem 1.9, p. 68), where the distance for two vector-functions $Y$ and $Z$ is given by
\[
d(Y, Z) = \sum\limits_{k=1}^m \max_{t\in [0,t_1]} ||y_k(t,\cdot)-z_k(t,\cdot)||_{L_2(\mathbb{R}^n)}.
\]
Let us denote the right-hand  side of equation (\ref{Volterra}) by $PU(t, x)$. To apply the Banach fixed point theorem we have to prove the following:

(a) if $U(t, x)\in {\mathbf C} ([0,t_1]; {\mathbf L_2}(\mathbb{R}^n))$, then $PU(t, x)\in {\mathbf C} ([0,t_1]; {\mathbf L_2}(\mathbb{R}^n))$;

(b) for any $Y, Z \in {\mathbf C} ([0,t_1]; {\mathbf L_2}(\mathbb{R}^n))$ one has
\[
d(P Y,P Z)\leq \delta\cdot d(Y, Z), \,\, \delta<1.
\]

Lemma \ref{WestimateH} implies condition (a). On the other hand, thanks to Lemma \ref{LestimateH} we arrive at
\[
    \bigg|\bigg|\int\limits_{0}^{t_1} S'(t_1-\eta, D)  \big[H(\eta,x, Y)-H(\eta,x, Z)\big] d\eta\bigg|\bigg|_{{\mathbf L}_2(\mathbb{R}^n)}\leq 
    \delta\cdot d(Y, Z),
    \]
where $\delta =c_1 t_1^\beta<1$ since condition (\ref{c1}).

Hence by the Banach fixed point theorem, there exists a unique solution  $U(t,x) \in  {\mathbf C} ([0,t_1]; {\mathbf L_2}(\mathbb{R}^n))$ to the integral equation (\ref{Volterra}). Next we consider the interval $[t_1, 2 t_1]$. Rewrite the equation (\ref{Volterra}) in the form
\[
 U(t,x) = F_1(t,x)+\int\limits_{t_1}^{t} S'(t-\eta, D) H(\eta,x, U) d\eta, \quad t_1\leq t\leq 2t_1, \ x\in \re^n.
\]
where
\[
F_1(x,t)=F(x,t)+ \int\limits_{0}^{t_1} S'(t-\eta, D) H(\eta,x, U) d\eta,
\]
is a known function, since the function $U(t,x)$ is uniquely defined on the interval $[0, t_1]$. Using the same arguments
as above, we derive that there exists a unique solution $U(t,x) \in  {\mathbf C} ([t_1, t_2]; {\mathbf L_2}(\mathbb{R}^n))$ to equation (\ref{Volterra}) on the interval $[t_1, 2t_1].$ Repeating this process, we conclude that there exists a unique solution $U(t,x) \in  {\mathbf C} ([0, T]; {\mathbf L_2}(\mathbb{R}^n))$ to equation (\ref{Volterra}).

We also need the following form of Gronwall's inequality (see, e.g., \cite{Kai}):
\begin{lemma}\label{GronLemma}
    Let $0<\rho<1$. Let the non-negative  function $y(t)\in C[0,T]$ and the positive constants $K_0$ and $K_1$ satisfy the inequality
   \begin{equation}\label{Gron1}
    y(t)\leq K_0+ \frac{K_1}{\Gamma(\rho)}\int\limits_0^t(t- s)^{\rho-1} y(s) ds
    \end{equation}
    for all $t\in [0, T]$. Then 
    \begin{equation}\label{Gron}
        y(t)\leq K_0\, E_\rho(K_1  t^\rho), \,\, t\in [0,T].
    \end{equation}

\end{lemma}
The appearance of the Mittag-Leffler function $E_\rho$ is explained by the fact that if in (\ref{Gron1}) inequality is replaced by equality, then the resulting integral equation is equivalent to the Cauchy problem $D_t^\rho   y(t)= K_1 y(t),\,\, y(0)= K_0$.

\begin{lemma}\label{ULtau} Let $F(t, x)= S(t,D) \varPhi(x)$ in equation (\ref{Volterra}) and $U(t, x)$ is the unique solution of this equation. Moreover, let the initial functions satisfy 
 the condition 
$
 \varPhi(x) \in {\mathbf L}^{\tau+\tau^*}_2(\mathbb{R}^n), \,  \tau > \frac{n}{2},$
and $H(t,x, U)$ satisfy conditions (\ref{H_L}) and (\ref{H0}).Then $U(t, x)\in {\mathbf C} ([0, T]; {\mathbf L}_2^{\tau+\tau^*}(\mathbb{R}^n)) $.
    
\end{lemma}
\emph{\bf Proof.} We have
\begin{equation}\label{IE}
    U(t,x) = S(t,D)\varPhi(x)+\int\limits_0^t S'(t-\eta, D) H(\eta,x, U) d\eta, \quad 0\leq t\leq T, \ x\in \re^n,
\end{equation}
or in component-wise form
\[
u_j(t,x)=v_j(t,x)+w_j(t,x, U),
\]
where
\[
v_j(t,x) = \frac{1}{(2\pi)^n} \sum_{k=1}^m  \sum_{q=1}^m \int\limits_{\re^n} e^{i\xi x}\mu_{j, q}(\xi)E_{\beta} \left( -\lambda_{q}(\xi)t^{\beta} \right) \nu_{q, k}(\xi) F[ \varphi_k] (\xi) d\xi,
\]
and
\[
w_j(t,x, U)=  \frac{1}{(2\pi)^n} \sum_{k=1}^m  \sum_{q=1}^m \int\limits_0^t \int\limits_{\mathbb{R}^n} e^{i\xi x}\mu_{j, q}(\xi)\eta^{\beta-1} E_{\beta, \beta} \left( -\lambda_{q}(\xi)\eta^{\beta} \right) \nu_{q, k}(\xi) F[ h_k] (t-\eta, \xi, U) d\xi d\eta. 
\]
As in (\ref{estimateKernel}), we see that $|\mu_{j, q}(\xi)E_{\beta} \left( -\lambda_{q}(\xi)t^{\beta} \right) \nu_{q, k}(\xi)|\leq C_{R_0, T}$.
By virtue of this estimate and (\ref{estimateKernel}), it is not difficult to prove (see Lemma \ref{WestimateH}) that
\[
||v_j(t,x)||_{L_2^{\tau+\tau^*}(\mathbb{R}^n)}\leq C\sum\limits_{k=1}^m ||\varphi_k||_{L_2^{\tau+\tau^*}(\mathbb{R}^n)},
\]
and 
\[
||w_j(t,x, U)||_{L_2^{\tau+\tau^*}(\mathbb{R}^n)}\leq C\max_{t\in [0,T]}\sum\limits_{k=1}^m ||h_k(t, \cdot, 0)||_{L_2^{\tau+\tau^*}(\mathbb{R}^n)}+ C\int\limits_0^t (t-\eta)^{\beta-1} \sum_{k=1}^m||u_k(\eta, \cdot)||_{L_2^{\tau+\tau^*}(\mathbb{R}^n)}.
\]

If we denote $y(t)=\sum_{k=1}^m||u_k(t, \cdot)||_{L_2^{\tau+\tau^*}(\mathbb{R}^n)}$, then by virtue of equation (\ref{IE}) we will have
\[
 y(t)\leq K_0+ \frac{K_1}{\Gamma(\beta)}\int\limits_0^t(t- s)^{\beta-1} y(s) ds, \,\, t\in [0, T],
    \]
where $K_1$ is a constant and 
\[
K_0=C\sum\limits_{k=1}^m ||\varphi_k||_{L_2^{\tau+\tau^*}(\mathbb{R}^n)}+C\max_{t\in [0,T]}\sum\limits_{k=1}^m ||h_k(t, \cdot, 0)||_{L_2^{\tau+\tau^*}(\mathbb{R}^n)}.
\]
Therefore, Gronwall’s inequality (\ref{Gron}) implies for any $j=1, \cdots, m,$
\begin{equation}\label{estimateU_j}
    ||u_j(t,x)||_{L_2^{\tau+\tau^*}(\mathbb{R}^n)}\leq  C\sum\limits_{k=1}^m ||\varphi_k||_{L_2^{\tau+\tau^*}(\mathbb{R}^n)}+C\max_{t\in [0,T]}\sum\limits_{k=1}^m ||h_k(t, \cdot, 0)||_{L_2^{\tau+\tau^*}(\mathbb{R}^n)}, \,\, t\in [0, T].   
\end{equation}
Lemma is proved.

\begin{theorem}
\label{Ntfph} Let matrix $\mathcal{A}(\xi)$ satisfy Conditions (A). Let the initial functions satisfy 
 the condition 
$
 \varPhi(x) \in {\mathbf L}^{\tau+\tau^*}_2(\mathbb{R}^n), \,  \tau > \frac{n}{2},$
and $H(t,x, U)$ satisfy conditions (\ref{H_L}) and (\ref{H0}).
Then the Nonlinear Problem has a unique solution $U(t,x)=\langle{u_1(t,x), u_2(t,x), \dots, u_m(t,x)}\rangle$, and each $u_j$  has the property
\[
\lim\limits_{|x|\rightarrow\infty} D^\alpha u_j(t, x)=0,\quad
|\alpha|\leq \ell_{j,j}, \quad t>0.
\]
Moreover, if $U^1=\langle{u_1^1, \cdots, u_m^1\rangle}$ and $U^2=\langle{u_1^2, \cdots, u_m^2\rangle}$ are two solutions to the Nonlinear Problem corresponding to two initial functions $\varPhi^1=\langle{\varphi_1^1, \cdots, \varphi_m^1\rangle}$ and $\varPhi^2=\langle{\varphi_1^2, \cdots, \varphi_m^2\rangle}$, respectively, then the stability estimate 
\begin{equation}\label{stability}
   \max_{t\in [0,T]}||u^1_{j}(t,x)-u_j^2 (t,x)||_{C(\mathbb{R}^n)}\leq C \sum\limits_{k=1}^m ||\varphi^1_k-\varphi_k^2||_{L_2^{{\tau+\tau^*}}(\mathbb{R}^n)},\,\, j=1, \cdots, m,
 \end{equation}
is valid.
\end{theorem}

\emph{\bf Proof.} Let us suppose that $U(t,x )$ is a solution of the Nonlinear Problem and $H(t,x, U)$ as a function of $(t,x)$ satisfies conditions
 \begin{equation}\label{condition_h}
   || h_j(t,x, U)||_{L_2^{\tau+\tau^*}(\mathbb{R}^n)} \in C[0,T], \,\, j=1, \cdots, m.
\end{equation}
  Then according to Theorem \ref{tfph}, $U(t,x )$ is a solution of the integral equation (see (\ref{Volterra}))
\[
 U(t,x) = F(t,x)+\int\limits_{0}^{t} S'(t-\eta, D) H(\eta,x, U) d\eta, \quad 0\leq t\leq T, \ x\in \re^n,
\]
where 
\[
F(t,x)=S(t,D)\varPhi(x).
\]
By estimate (\ref{estimateV_0}) one has $F(t,x)\in {\mathbf C} ([0,T]; {\mathbf L_2}(\mathbb{R}^n))$. Therefore, this integral equation, by virtue of the Lemma  \ref{VolterraSolution}, has a unique solution $U(t,x) \in  {\mathbf C} ([0,T]; {\mathbf L_2}(\mathbb{R}^n))$.

It remains to prove that the unique solution to equation (\ref{Volterra}) that we have found satisfies all the conditions of the Nonlinear Problem and $H(t,x, U)$ satisfies condition (\ref{condition_h}).

As for (\ref{condition_h}), we have by (\ref{H_L})-(\ref{0_H_L})
\[
|| h_j(t,x, U)||_{L_2^{\tau+\tau^*}(\mathbb{R}^n)}\leq || _0 h_j(t,x, U)||_{L_2^{\tau+\tau^*}(\mathbb{R}^n)}+|| h_j(t,x, 0)||_{L_2^{\tau+\tau^*}(\mathbb{R}^n)}\leq
\]
\[
L_0 ||u_j(t,x)||_{L_2^{\tau+\tau^*}(\mathbb{R}^n)}+|| h_j(t,x, 0)||_{L_2^{\tau+\tau^*}(\mathbb{R}^n)}\leq
\]
(and by (\ref{estimateU_j}))
\[
C\sum\limits_{k=1}^m ||\varphi_k||_{L_2^{\tau+\tau^*}(\mathbb{R}^n)}+C\max_{t\in [0,T]}\sum\limits_{k=1}^m ||h_k(t, \cdot, 0)||_{L_2^{\tau+\tau^*}(\mathbb{R}^n)}, \,\, t\in [0, T].
\]

Using this estimate, all the conditions of the Nonlinear Problem are proved in the same way as in Theorem \ref{tfph}.

The stability estimate (\ref{stability}) is a consequence of the Lipschitz condition for function $H$ and Lemma \ref{GronLemma}.

\begin{remark}\label{nonlinearNessesary}
When $\tau^*=0$, a sufficient condition (\ref{H0}) on the function $H(t,x, 0)$ is also necessary for the existence and uniqueness of a solution to Nonlinear Problem for the same reason as in Remark \ref{nessesary}.
\end{remark}
\begin{remark}\label{nonlinear}To our best
knowledge Nonlinear Problem for  a system of fractional order partial differential equation has been considered for the first time.
\end{remark}

\section{Conclusions}
The paper considers a system of fractional partial differential equations in the case when the matrix of the system $\mathcal{A}(\xi)$  satisfies Condition $(A)$ and the order of fractional derivatives is a scalar $\beta\in (0,1]$, i.e. same for all equations of the system. Sufficient conditions on the initial function and on the right-hand side of the equation are found that guarantee the existence and uniqueness of classical solutions for Linear and Nonlinear Problems. In some cases these conditions are also necessary (see Remarks \ref{nessesary} and \ref{nonlinearNessesary}). 

Let us list what generalizations of these results can be expected: 1) without much difficulty one can consider “hyperbolic” systems, i.e. $\beta \in (1,2]$, adjusting properly the initial conditions; 2) these results can also be easily transferred to the case of fractional derivatives in the Riemann-Liouville sense; 3) we can reject the Hermitian nature of matrix $\mathcal{A}(\xi)$  (see, for example, \cite{Varsha},  \cite{Odabat} and \cite{UmarovToAppear}); 4) we can consider in the system (\ref{system_main}) a vector of order $\beta=\langle{\beta_1, \cdots, \beta_m\rangle}$ of fractional derivatives (i.e., each equation of the system has its own order of fractional derivatives), where $\beta_j$ are rational numbers. It is important to note that, as can be seen from the results of work \cite{UmarovToAppear}, after reducing the original system to a system with the same order of fractional derivatives, the eigenvalues of the original matrix participate in the formula for representing the solution and the same asymptotic estimates (\ref{lambdaEstimate}) are valid for them.

These generalizations are the subject of further research.


\begin{thebibliography}{99}

\bibitem{DasGupta} S. Das, P.K. Gupta,  A mathematical model on fractional Lotka-Volterra equations. Journal of theoretical biology, {\bf 277} (1), 1-6, 2011.



\bibitem{Rihan} F.Rihan,
Numerical Modeling of Fractional-Order Biological Systems. Abstract and Applied Analysis, {\bf 2013}, 1--13, 2013.





\bibitem{GuoFang} Ch.Guo, Sh. Fang, Stability and approximate analytic solutions of the fractional Lotka-Volterra equations for three competitors. Advanced difference equations, 219, 1-14, 2016.






\bibitem{Khan} N.A. Khan, O.A. Razzaq, S.P. Mondal, Q. Rubbab,
Fractional order ecological system for complexities of interacting species with harvesting threshold in imprecise environment. Advances in Difference Equations, 405, 1-34, 2019.





\bibitem{Rana} S. Rana, S. Bhattacharya, J. Pal, Guerekata G., Chattopadhyay.
Paradox of enrichment: A fractional differential approach with memory. Physica A: Statistical Mechanics and its Applications, {\bf 392} (17), 3610--3621, 2013.





\bibitem{Zeb} A. Zeb, G. Zaman, M.I. Chohan, Sh. Momani, V.S. Erturk,
Analytic numeric solution for SIRC epidemic model in fractional order.  Asian J. of Math and Appl. {\bf 2013}, 1-19, 2013.





\bibitem{Islam} R. Islam, A. Pease, D. Medina, Oraby T. 
Integer Versus Fractional Order SEIR Deterministic and Stochastic Models of Measles. International Journal of Environmental Research and Public Health, {\bf 17} (6), 1-19, 2020.




\bibitem{Lyapunov-2} Vargas-de-Le\'on. Volterra-type Lyapunov functions for fractional-order epidemic systems. Coomun. Nonlinear Sci. Numer. Simulat., {\bf 24}, 75--85, 2015.





\bibitem {Goufo} E.F.D. Goufo, R. Maritz and J. Munganga, Some properties of Kermack-McKendrick epidemic model with fractional derivative and nonlinear incidence. Adv. Difference Equ. {\bf 2014} (1), Article ID 278, 1--9, 2014.



\bibitem{Almeida} R. Almeida, Analysis of fractional SEIR model with treatment. Applied Mathematical Letters, {\bf 84} 56--62, 2018.





\bibitem{Rajagopal} K. Rajagopal, N. Hasanzadeh, F. Parastesh, I.I. Hamarash, S. Jafari, I. Hussain,  A fractional-order model for the novel coronavirus (COVID-19) outbreak. Nonlinear Dynamics, {\bf 101}, 701--718, 2020.






\bibitem{Varsha} V. Daftardar-Gejji, A. Babakhani Analysis of a system of fractional differential equations. J. of Math. Anal. and Appl., {\bf 293}, 511-522, 2004.




\bibitem{BonillaKilbas} B. Bonilla, A.A. KIlbas, J.J. Trujillo, Systems of nonlinesr fractional differential equations in the space of summable functions, Proceedings of the Institute of Mathematics of the National Academy of Sciences of Belarus, 2000, V. 6, 38-46.




\bibitem{DengLiGuo} W. Deng, Ch. Li, Q. Guo, Analysis of fractional differential equations with multi-orders.  Fractals, {\bf 15} (2), 173--182, 2007.





\bibitem{Odabat} Z. Odibat, Analytic study on linear systems of fractional differential equations. Computers and Mathematics with Applications. {\bf 59}, 1171--1183, 2010.





\bibitem{Mamchuev} M.O. Mamchuev, Boundary value problem for a multidinensional system of
equations with Riemann–Liouvile fractional derivatives, Sib. Elektron. Mat. Izv. ` ,
2019, Volume 16, 732–747




\bibitem{NewMethod} H. Jafaria, M. Nazarib, D. Baleanuc, C.M. Khalique, A new approach for solving a system of fractional partial
differential equations, Computers and Mathematics with Applications, 66 (2013), 838-843.





\bibitem{Erturk}
V.S. Ert\"urk,  S. Momani, Solving systems of fractional differential equations using differential transform method. J. Comput. Appl. Math., {\bf 215}, 142-151, 2008. 









\bibitem{Kochubei2} N. Anatoly Kochubei,
Fractional-hyperbolic equations and systems. Cauchy problem, In ”Handbook of Fractional Calculus”,  V. 2, pp. 197-223, De-Gruyter, 2019.













\bibitem{Abdulaziz} O. Abdulaziz, I. Hashim, S. Momani, Solving systems of fractional differential equations by homotopy-perturbation method. Physics Letters A, {\bf 372}, 451--459, 2008.






\bibitem{Lyapunov-1} N. Aguila-Comacho, Duarte-Mermoud M.A., Gallegos J.A. Lyapunov functions for fractional order systems. Commun. Nonlinear. Sci. Simulat. {\bf 19}, 2951--2957, 2014.


\bibitem{Wang} G. Wang, R.P. Agarwal, A. Cabada, Existence results and the monotone iterative technique for systems of nonlinear
fractional differential equations. Applied Mathematics Letters, {\bf 25} (6), 1019--1024, 2012.




\bibitem{UAChen} S.R. Umarov, R.R. Ashurov, YangQuan Chen. On a method of solution of systems of fractional pseudo-differential equations. // Fractional Calculus and Applied Analysis, 2021, V. 24, № 1, pp. 254-277.





\bibitem{UmarovToAppear} S. Umarov,  Representations of solutions of systems of time-fractional pseudodifferential equations, Frac. Calc. Appl. Anal. (2024) (to appear)







\bibitem{UmarovNew} S. Umarov, Representations of solutions of
time-fractional multi-order systems of
differential-operator equations, arXiv:submit/ 5388903 [math-ph] 4Feb 2024


\bibitem{Kochuei1}N. Anatoly Kochubei
Fractional-parabolic equations and systems.
Cauchy problem, In ”Handbook of Fractional Calculus”, V. 2, pp. 145-159, De-Gruyter, 2019.




\bibitem{Kochubei3} A.N. Kochubei, Fractional-parabolic systems, Potential Anal., 37 (2012), 1–30.





\bibitem{Kochubei4} A.N. Kochubei, Cauchy problem for fractional diffusion-wave equation with variable
coefficients, Appl. Anal., 93 (2014), 2211–2242.

\bibitem{Kochubei2} N. Anatoly Kochubei,
Fractional-hyperbolic equations and systems. Cauchy problem, In ”Handbook of Fractional Calculus”,  V. 2, pp. 197-223, De-Gruyter, 2019.






\bibitem{Eidelman}S.D. Eidelman, Parabolic Systems, North-Holland, Amsterdam, 1969.





\bibitem{Krein} S.G. Krein, Linear Differential Equations in Banach Space, American Mathematical Society,
Providence, 1972.






\bibitem{GelShil}I.M. Gelfand, G.E. Shilov, Generalized functions, V. 3, Theory of Differential Equations, Academic Press, New York, San Francisco, London, 1967. 






\bibitem{Gantmacher}F.R. Gantmacher: The theory of matrices, V. 1, Chelsea Publishing Company, New York, 1959.





\bibitem{Bez} E.G. Bazhlekova, Subordination principle for fractional evolution equations, Fract. Calc. Appl.
Anal., 3 (2000), 213–230.





\bibitem{KST} A.A. Kilbas, H.M. Srivastava, 
J.J. Trijillo,  Theory and Applications of Fractional Differential Equations. Elsevier Science, 2006.





\bibitem{Horn} R.A. Horn, C.R. Johnson. Matrix Analysis, Cambridge University Press, Ch. 7, 1985.







\bibitem{Bhatia} R. Bhatia, Positive definite matrices, Princeton Series in Applied Mathematics, 2007.














\bibitem{Handbook} A. Kochubey, Yu. Luchko, De Gruyter, Handbook of Fractional Calculus with Applications. Volume 2: Fractional Differential Equations.  2019.







\bibitem{Petrowsky}I. Petrowsky, Über das Cauchysche Problem fur ein System linearer partieller
Differentialgleichungen im Gebiete der nichtanalytischen Funktionen,Bull. Univ. Etat. Moscou,
Ser. Int., Sect. A, Math. et Mécan. 1, Fasc.,7(1938), 1–74.

\bibitem{Bellman} R. Bellman, Matrix analysis, SIAM, Philadelphia, USA.






\bibitem{Voevodin} V.V. Voevodin, Computational foundations of linear algebra, Nauka, Moscow, 1977.













\bibitem{Duhamel_1} S. Umarov, On fractional Duhamel's principle and its applications. Journal of Differential Equations, {\bf 251} (10), 5217--5234, 2012.






\bibitem{Umarov_book_2015} S. Umarov, Introduction to Fractional and Pseudo-Differential Equations with Singular Symbols. Springer, 2015.






\bibitem{AZRn} R.R. Ashurov, R.T. Zunnunov, Initial-boundary value and inverse problems for subdiffusion equation in
$R^n$. Fractional differential calculus, {\bf 2020}.{\em 10},.291-306.






\bibitem{Kai} K. Diethelm, The Analysis of Differential Equations of Fractional Order: An
Application-Oriented Exposition Using Differential Operators of Caputo Type.
Lecture Notes in Mathematics, vol. 2004. Springer, Berlin/Heidelberg (2010)












    










\bibitem{Baz} E.G. Bazhlekova,  Subordination principle for fractional evolution equations. Fract. Calc. Appl.
Anal. 3, 213–230 (2000)

\bibitem{Podlubny} I. Podlubny,  Fractional Differential Equations. Academic Press, 1998.



\bibitem{SKM} S.G. Samko, A.A. Kilbas, O.I. Marichev,   Fractional Integrals and Derivatives: Theory  and Applications.  Gordon and Breach Science Publishers, 1993.

\end{thebibliography}
\end{document}